\documentclass[12pt,english]{article}
\usepackage[margin=3cm]{geometry}
\usepackage{subcaption}
\usepackage{amsfonts}

\usepackage{amsmath}
\usepackage{amsthm}
\usepackage{color}
\usepackage{graphicx}
\usepackage{dsfont}
\usepackage{soul}
\usepackage[ruled,vlined,onelanguage]{algorithm2e}
\usepackage{tikz}
\usepackage{array}
\usepackage{times}
\usepackage{xspace}

\usepackage{microtype} % added by SLD

\usepackage{makecell}
\usepackage{multirow}
\usepackage{booktabs}
\usepackage{rotating}

\usepackage{float}

\usepackage{hyperref}
\hypersetup{
    unicode      = false,
    pdftoolbar   = true,
    pdfmenubar   = true,
    pdffitwindow = true,
    pdftitle     = {Uncertainty with ensembles of surrogates},
    pdfauthor    = {Audet et al.},
    pdfsubject   = {ensembles of surrogates},
    pdfkeywords  = {ensembles of surrogates},
    pdfnewwindow = true,
    colorlinks   = true,
    linkcolor    = blue,
    citecolor    = blue,
    filecolor    = black,
    urlcolor     = blue,
    breaklinks   = true
}
\newcommand{\mads}{{MADS}\xspace}
\newcommand{\orthomads}{{OrthoMADS}\xspace}
\newcommand{\nomad}{{\sf NOMAD}\xspace}
\newcommand{\solar}{{\sf solar}\xspace}
\newcommand{\styrene}{{\sf styrene}\xspace}
\newcommand{\dfn}{{\sf DFN}\xspace}
\newcommand{\shebo}{{\sf SHEBO}\xspace}

\definecolor{Red}{rgb}{1,0,0}
\definecolor{Green}{rgb}{0,.6,0}
\definecolor{Blue}{rgb}{0,0,1}

\title{Quantifying uncertainty with ensembles of surrogates for blackbox optimization
\thanks{This work is supported by the IVADO Fundamental Research Project Grant PRF-2019-8079623546.}}
\author{
    \href{mailto:Charles.Audet@gerad.ca}{Charles~Audet}
    \thanks{
       \href{https://www.gerad.ca}{GERAD} and \href{https://www.polymtl.ca}{Polytechnique Montr\'eal},
          \href{https://www.gerad.ca/Charles.Audet}{\tt www.gerad.ca/Charles.Audet}
  }
  \and
        \href{mailto:Sebastien.Le.Digabel@gerad.ca}{S\'ebastien~{Le~Digabel}}
        \thanks{
        \href{https://www.gerad.ca}{GERAD} and \href{https://www.polymtl.ca}{Polytechnique Montr\'eal},
                  \href{https://www.gerad.ca/Sebastien.Le.Digabel}{\tt www.gerad.ca/Sebastien.Le.Digabel}
  }
  \and
  \href{mailto:renaud.saltet@polymtl.ca}{Renaud~Saltet}
  \thanks{
  \href{https://www.gerad.ca}{GERAD} and \href{https://www.polymtl.ca}{Polytechnique Montr\'eal}, \href{mailto:renaud.saltet@polymtl.ca}{\tt renaud.saltet@polymtl.ca}
  }
}
\begin{document}

\maketitle
\thispagestyle{empty}

%-------------------------------------------------------------%
\noindent
{\bf Abstract:}
{
This work is in the context of blackbox optimization where the functions defining the problem are expensive to evaluate and where no derivatives are available. A tried and tested technique is to build surrogates of the objective and the constraints in order to conduct the optimization at a cheaper computational cost. This work proposes different uncertainty measures when using ensembles of surrogates. The resulting combination of an ensemble of surrogates with our measures behaves as a stochastic model and allows the use of efficient Bayesian optimization tools. The method is incorporated in the search step of the mesh adaptive direct search (\mads) algorithm to improve the exploration of the search space. Computational experiments are conducted on seven analytical problems, two multi-disciplinary optimization problems and two simulation problems. The results show that the proposed approach solves expensive simulation-based problems at a greater precision and with a lower computational effort than stochastic models.
} \\

%-------------------------------------------------------------%

%-------------------------------------------%
\section{Introduction}
%-------------------------------------------%

This work considers the constrained optimization problem
\begin{equation}\label{eq:P}
\tag{$P$}
    \begin{aligned}
    \min_{x\in\mathcal{X}}&\ \ f(x)\\
    \mbox{s.t.}&\ \ x\in\Omega
    \end{aligned}
\end{equation}
where $\mathcal{X}$ is a subset of $\mathbb{R}^n$; $\Omega$ denotes the feasible set $\big\{x\in\mathcal{X}\ |\ c_j(x)\leq0,\ j\in\{1,2,\dots,m\}\big\}$, where the functions $c_j:\mathbb{R}^n \rightarrow \overline{\mathbb{R}} = \mathbb{R}\cup\{-\infty,+\infty\}$ are the constraint functions of the problem; and $f:\mathbb{R}^n \rightarrow \overline{\mathbb{R}}$ is the objective function. The set $\mathcal{X}$ contains the points that satisfy \textit{unrelaxable} constraints~\cite{LedWild2015}: every point explored during the optimization process must lie in $\mathcal{X}$ either because $f$ is not defined elsewhere or because a point outside $\mathcal{X}$ has no meaning in the original problem, e.g., a negative length or a probability greater than one. The set $\mathcal{X}$ typically represents bound constraints of the form $\mathcal{X}=\{x\in\mathbb{R}^n\ | \ \ell\leq x\leq u\}$ where $\ell$ and $u$ are vectors of $\{\mathbb{R}\cup\{-\infty\}\}^n$ and $\{\mathbb{R}\cup\{+\infty\}\}^n$, respectively. The functions $c_j$,~$j\in\{1,2,\dots,m\}$, denote \textit{relaxable} constraints, which means that they can be violated during the optimization process, however, the final solution must satisfy these constraints.

In \textit{blackbox optimization} (BBO), no information is available on the functions $f$ and $c_j$,\linebreak $j~\in~\{1,2,\dots,m\}$, beyond the mere values they produce, hence the \textit{blackbox} designation. In particular, no derivatives can be used either because they are especially hard to estimate or because they do not exist. Designing algorithms that do not use derivatives is referred to as \textit{Derivative-Free Optimization} (DFO). This context typically occurs when the functions are the results of numerical simulations. Consequently, a blackbox is assumed to be costly, i.e., one evaluation might take seconds~\cite{ChEsKeDaEt00a}, minutes~\cite{AuBeCh2008a,MuPaSaVaArFaAg2020}, hours~\cite{AuOr06a,Brass2017} or even days~\cite{MaWaDeMo07}. BBO consists in designing algorithms capable of finding the best possible solution to such a problem with a given budget of function evaluations. For a better understanding of the theoretical importance of the existence of derivatives in DFO, see~\cite{AuHa2020}. Two reference books are available in DFO and BBO~\cite{AuHa2017,CoScVibook} as well as an extended review~\cite{LaMeWi2019}.

BBO algorithms can be roughly divided in two categories: \textit{direct-search} methods and methods using \textit{surrogates}. Direct-search algorithms only use comparisons between points and no other information like an approximation of derivatives. This philosophy has led to fruitful algorithmic frameworks such as the state-of-the-art algorithms \textit{generalized pattern search} (GPS)~\cite{Torc97a} and \textit{mesh-adaptive direct search} (\mads)~\cite{AuDe2006}, the latter is described in this article. Both frameworks lie on the search-poll paradigm~\cite{Audet04a}: the \textit{search} step offers flexibility for the user to implement any method they see fit for the problem to optimize, while the \textit{poll} step imposes more rigid procedures in order to guaranty convergence and explore further the surroundings of the best known solution. The second category resorts to \textit{surrogates}, i.e., functions that are expected to mimic the behaviour of the objective and the constraints while being significantly cheaper to evaluated. A surrogate can be a simplified and static version of the blackbox that do not evolve over the optimization, or a dynamic surrogate that is based on regression or interpolation on the previously evaluated points. A dynamic surrogate is also called a model. One can then expect that minimizing a surrogate of the objective while satisfying surrogates of the constraints will lead to a promising new candidate point for the true problem. The combination of direct-search and model-based methods has proven fruitful. Surrogates can be incorporated either in a subproblem embedded in the search step to find candidates points, or as means of ranking candidate points of the poll step when opportunistic strategies are used.

Among the model-based methods, ensembles of models and stochastic models are two efficient techniques. Ensembles of models consist in giving each model a weight that is supposed to reflect its quality, and the combination of all weighted models yields an \textit{aggregate model} that can be used as a standard surrogate. Stochastic models not only produce a prediction at a given point but also a measure of the uncertainty on this prediction, which is fit for Bayesian optimization. Ensembles of models and stochastic models have both proven efficient but remain fundamentally separated. Using several deterministic surrogates naturally produces a deterministic aggregate model which is incompatible with Bayesian optimization. If some of the models used are stochastic though, the provided uncertainty can be exploited by all the models as in~\cite{ViHaWa2013}. However, typical stochastic surrogates as Gaussian processes (GPs) become particularly costly to train as the training set grows. The proposed approach is to identify areas where the predictions of the models differ from each other in order to derive some form of uncertainty. The idea of using the correlation between several models to guide the optimization has already been tackled in the literature. In~\cite{GoTuHaShQu2007}, the deviation between the predicted values of several models is used \textit{a posteriori} to check the overall quality of the aggregate model. In~\cite{Mu2020}, two surrogates of the objective are available: a low-fidelity one that is cheap to compute but not accurate and a high-fidelity one that is in contrast more expensive and more reliable. Then RBF models of the two surrogates are computed and the correlation between them is used to choose which surrogate, high or low-fidelity, to evaluate next. In~\cite{To2015}, the correlation between variable fidelity co-kriging models is exploited to determine which model to evaluate. In~\cite{RuJiZhHuSh2020}, the correlation between variable fidelity multi-level generalized co-kriging models is incorporated in an extended probability of improvement.

The contribution of this work is an extension to ensembles of models when used in the form of aggregate models. For a point $x$ of the search space $\mathcal{X}$, the extended aggregate models produce not only a prediction $\hat f(x)$, but also an uncertainty $\hat\sigma(x)$, therefore imitating a stochastic model. The resulting surrogate is then exploited in the search step of \mads in subproblems inspired by Bayesian optimization. The proposed approach has been tested on seven analytical problems, two multi-disciplinary optimization problems and two simulation problems. It has been compared to other versions of \mads as well as two other BBO solvers. Results show that the proposed extended aggregate models manage to find solutions of most of the difficult real-world problems at a greater precision than the other algorithms, and with less computational effort than the competing stochastic models.

The manuscript is structured as follows. Section~\ref{sec:background} introduces ensemble of surrogates and the Bayesian optimization framework as well as a high-level description of the \mads algorithm. Section~\ref{sec:contribution} describes the quantification of uncertainty when using ensembles of models and our incorporation of the resulting extended aggregtae models into the \mads algorithm. Section~\ref{sec:results} shows the computational results on the set of problems. A concluding discussion is proposed in Section~\ref{sec:conclusion}.

%-------------------------------------------%
\section{Background}
\label{sec:background}
%-------------------------------------------%

This section describes the use of surrogates in BBO, with a special focus on stochastic surrogates and ensembles of surrogates, as well as the \mads algorithm.

\subsection{Surrogates in BBO}\label{subsec:surrogates}

A common approach in BBO uses \textit{surrogates} of the objective and the constraints in order to guide the optimization. A surrogate shares similarities with the true functions of the problem while being significantly cheaper. Two types of surrogates can be distinguished: \textit{static surrogates} and \textit{dynamic models}. A static surrogate is a simplified version of the blackbox that can be obtained for example through a simplified physics model, a coarser mesh in a finite elements simulation or a looser stopping criterion in a numerical method. Such a surrogate is fixed and does not evolve over the optimization, hence the \textit{static} designation. On the other hand, a dynamic model is an interpolation or regression model that approaches the true functions by fitting previously evaluated sample points. Since it attempts to approximate the true function, the \textit{model} designation is more appropriate than surrogate in this case. Common models used in BBO are polynomial response surfaces (PRS)~\cite{AcarOptimizedWeight2009,MuPi2011} and especially quadratic models~\cite{CoLed2011,ARConn_PhLToint_1996,CuRoVi10}, radial basis functions (RBF)~\cite{MBjorkman_and_KHolmstrong_2000,KiArYa2011,Regis2011Stoch,Regis2014HighDim,RGRegis_CAShoemaker_2005,AnVaVi2009,WiReSh2008,WiSh2011}, support vector machines (SVM)~\cite{VeKa2017}, kernel smoothing (KS)~\cite{AcarOptimizedWeight2009}, and Gaussian processes (GPs)~\cite{Book00a,PeHiLaCo14,RaWi06,SaWeMiWy89a}.

Surrogates are basically used for two purposes in BBO: finding new candidate points and ranking existing candidate points before evaluation by the true problem. The surrogate management framework~\cite{BoDeFrSeToTr99a} establishes the interplay between the surrogate evaluations and the true evaluations and is described in Algorithm~\ref{algo:SMF} as in~\cite[Chapter~13]{AuHa2017}.
%
%---------------------------------------%
\begin{algorithm}[htbp]
\caption{Surrogate management framework.} 
\SetAlgoNoLine
\SetKwProg{subalg}{}{}{}

\subalg{\textbf{1. Exploration using the surrogate}}
{\SetAlgoVlined
Use the surrogate problem to generate a list $\mathcal{L}$ of candidate points\\
Evaluate the true functions at points in $\mathcal{L}$ in an opportunistic way\\
If a new incumbent solution is found, go to 3; otherwise go to 2
}

\subalg{\textbf{2. Ranking using the surrogate}}
{\SetAlgoVlined
Use the optimization algorithm to generate a list $\mathcal{L}$ of candidate points\\
Use the surrogate functions to order the points in $\mathcal{L}$\\
Evaluate the true functions at points in $\mathcal{L}$ in an opportunistic way\\
}

\subalg{\textbf{3. Parameters update}}
{\SetAlgoVlined
Update algorithmic parameters\\
Check stopping criteria or go to 4
}

\subalg{\textbf{4. Model update (optional)}}
{\SetAlgoVlined
Update the model by using the new values of the true functions obtained in 1 and 2
}
\label{algo:SMF}
\end{algorithm}
%---------------------------------------%

For extended reviews on model-based or model-assisted optimization, see~\cite{AuHa2020,VuDaHaLi2016} and~\cite[Section~2.2]{LaMeWi2019}.

%---------------------------------------%
\subsubsection{Ensemble of models}
\label{subsubsection:ensembles}
%---------------------------------------%

In an optimization context, there is not one type of models that dominates the others~\cite{AcarOptimizedWeight2009,GoTuHaShQu2007}. Even for a given problem, the best performance can be obtained with different models depending on the initial sampling~\cite{GoTuHaShQu2007}. A tempting strategy is to resort to several models simultaneously. This idea has proven efficient in several works~\cite{AcarOptimizedWeight2009,AuKoLedTa2016,Chen2018,GoTuHaShQu2007,MuPi2011,VianaEnsemble2008,YePanDong2018} in which an ensemble of models is used to build an aggregate model defined by
\begin{equation}\label{eq:aggregated_model}
    \hat f(x) = \sum_{p=1}^s w^p\tilde f^p(x)
\end{equation}
where $s\geq1$ is the number of models; $\tilde f^1,\tilde f^2,\dots,\tilde f^s$ are $s$ models of the objective $f$ built from sample points; and $w^1,w^2,\dots,w^s$ are positive weights such that $\sum_{p=1}^s w^p=1$. The main difficulty lies in the attribution of the weights that must reflect the quality of the models. To do so, weights can be attributed so that the error of the aggregate model will be minimal, thus introducing an optimization subproblem~\cite{AcarOptimizedWeight2009,VianaEnsemble2008}. Another approach is to compute an error metric $\mathcal{E}^p$ for each model $\tilde f^p$ and then attribute a weight $w^p$ that is a function of $\mathcal{E}^p$~\cite{AuKoLedTa2016,GoTuHaShQu2007,MuPi2011}. In this type of strategy, an error metric must be chosen first. Common metrics for this purpose are statistical measures with cross-validation like root mean square error (RSME) and predicted residual sum of squares (PRESS) that take into account the gaps between the values of the models and those of the true objective. In~\cite{AuKoLedTa2016}, the authors propose an error metric that is not statistical but is rather a measure of a model's capacity to rank points in the same order as the objective would do. The rationale for this latter metric is that in a BBO context a good model does not necessarily approximate well the values of the true function but is rather capable of discriminating candidate points and telling apart promising ones.

Once an error metric is chosen, weights must be attributed accordingly. For instance, in~\cite{GoTuHaShQu2007}, the authors propose the three following options:
$$
    \begin{aligned}
        &w^p \propto \mathcal{E}^{\mathrm{tot}}-\mathcal{E}^p\\
        \mbox{or } &w^p \propto \mathds{1}_{\mathcal{E}^p=\mathcal{E}^{\mathrm{min}}}\\
        \mbox{or } &w^p \propto (\mathcal{E}^p + \alpha\mathcal{E}^{\mathrm{av}})^{\beta}
    \end{aligned}
$$
where $\mathcal{E}^{\mathrm{tot}}$ is the total error of all models, $\mathcal{E}^{\mathrm{min}}$ is the minimal error, $\mathcal{E}^{\mathrm{av}}$ is the average error, and $\alpha<1$ and $\beta<0$ are adjustable parameters.

%-------------------------------------%
\subsubsection{Bayesian optimization}
%-------------------------------------%

In general, models may be trusted in areas where the true functions have been sufficiently sampled. But the further away from these areas, the less accurate are the models. In Bayesian optimization, the objective function is interpreted as a stochastic process: an \textit{a priori} distribution $P[f]$ is assumed, then with the set of sample points $\mathbb{X}$ and a likelihood model $P[\mathbb{X}\ |\ f]$, an \textit{a posteriori} distribution $P[f\ |\ \mathbb{X}]$ is built thanks to Bayes' rule. Consequently, for any point $x$ of the search space, a stochastic model not only produces a prediction $\mu(x)$ but also a measure of uncertainty on that prediction $\sigma(x)$. If properly exploited, this uncertainty enables to explore areas in which the model confesses to be unreliable, instead of spending the entire budget on restricted areas. This is called the compromise between exploration and exploitation. Commonly used stochastic models are generalized linear models~\cite{NeWe1972}, Gaussian processes~\cite{RaWi06}, and dynamic trees~\cite{TaGrPo2011}.

The compromise between exploration and exploitation is then realized with an acquisition function. A simple example is the upper confidence bound (UCB)~\cite{SrKrKaSe2009} that takes into account the most optimistic value, i.e., minimizes $\mu(x) - \kappa\sigma(x)$. A more sophisticated instance is the probability of improvement (PI)~\cite{DRJones_2001} which is the probability that the objective decreases from the best known value at a given point. Finally, a very popular example is the expected improvement (EI)~\cite{MoTiZi1978} that not only takes into account the probability of decrease but also the expected amplitude thereof.

\begin{figure}[htbp]
  \centering
  \includegraphics[width=0.8\linewidth]{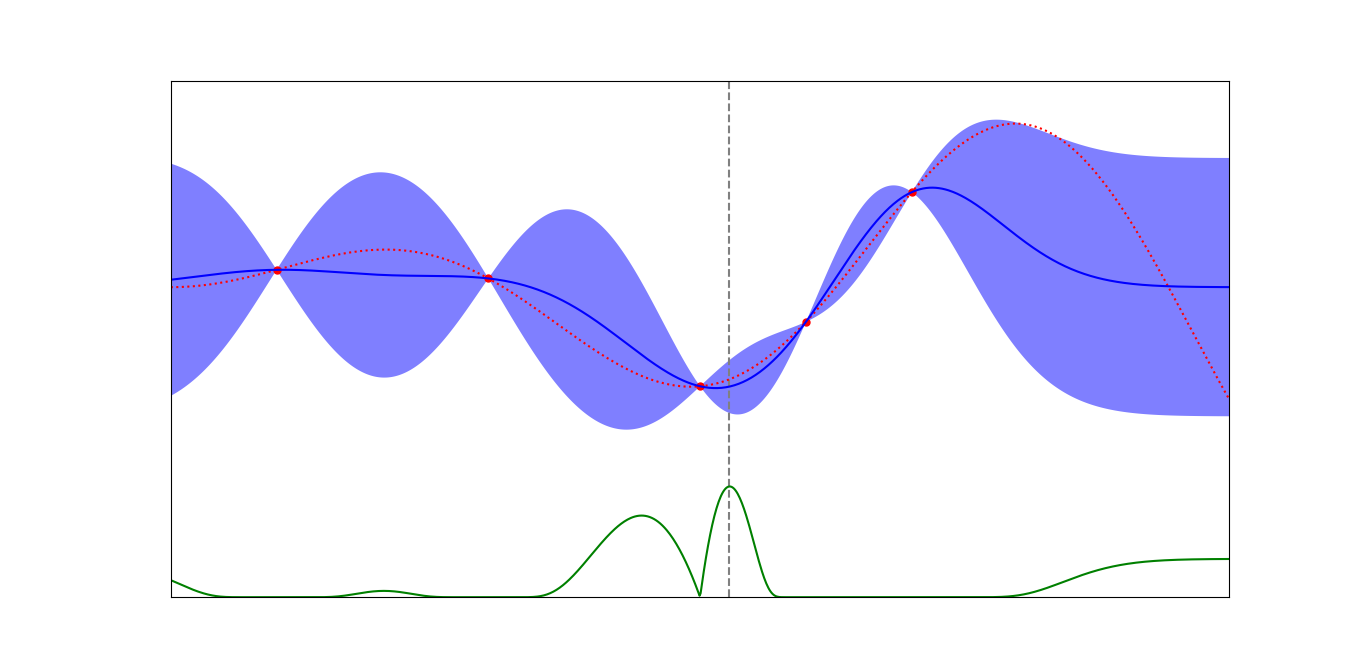}
  \setlength{\unitlength}{0.8cm}
  \begin{picture}(0,0)
  \put(-2.1,0.35){$x$}
  \put(-15.1,4.3){{\color[rgb]{0,0,0.7}$f(x)$}}
  \put(-1.6,1.1){{\color[rgb]{0.2,0.6,0.2}$EI(x)$}}
  \end{picture}
  \caption{Kriging and expected improvement (EI) on $f:x\mapsto x\sin x$.}
  \label{fig:GP_EI}
\end{figure}
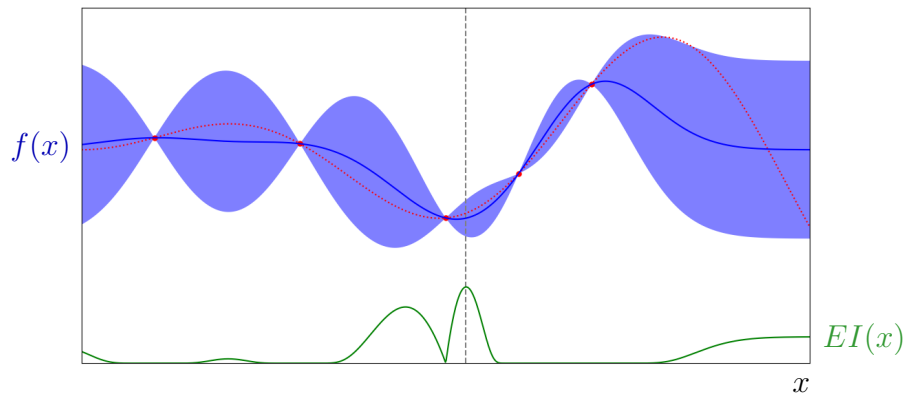

Figure~\ref{fig:GP_EI} shows an example of Gaussian process regression - also know as {\em kriging}~-  as well as the resulting expected improvement on a one-dimensional objective function. The dashed curve represents the objective $f:x\mapsto x\sin x$; the five dots are the sample points; the curve interpolating the dots is the prediction $\mu:x\mapsto\mu(x)$; an the filled area represents the 95\% confidence interval given at any point $x$ by $[\mu(x)+1.96\sigma(x),\ \mu(x)-1.96\sigma(x)]$. The curve at the bottom represents - in a different scale - the expected improvement EI. The resulting candidate point maximizes EI and is indicated by the vertical dashed line.
%
%--------------------------%
%--------------------------%
%
Maximizing EI is a method introduced by Jones et al.~\cite{JoScWe1998} that is efficient and easy to grasp. More elaborated techniques have since been developed. Talgorn et al.~\cite{TaLeDKo2014} propose various formulations for the surrogate subproblem that use multiple acquisition functions at once.

For an extended literature review on Bayesian optimization, see~\cite{GrRaGuVeVe2020,SSWAF2015}.

%-------------------------------------------%
\subsection[The MADS algorithm]{The \mads algorithm}
\label{subsec:mads}
%-------------------------------------------%

\mads~\cite{AuDe2006} is a direct-search algorithmic framework that follows a search-poll paradigm in which the mandatory poll step guaranties the convergence and the optional search step gives room for flexible exploration techniques. In order to ensure convergence, every candidate point must lie on a \textit{mesh} defined at Iteration~$k$ by
%
%-------------------------------------------%
$$M^k=\left\{ x+\delta^k Dy\ :\ x\in V^k,\  y\in\mathbb{N}^{n_D} \right\} \subset\mathbb{R}^n$$
%-------------------------------------------%
%
where $V^k\subset\mathbb{R}^n$ is the \textit{cache}, i.e., the set of all evaluated points up to Iteration $k$; $\delta^k>0$ is the \textit{mesh size parameter}; and $D$ is a fixed matrix of $\mathbb{R}^{n\times n_D}$, the columns of which represent $n_D$ directions of $\mathbb{R}^n$. Before the algorithm starts, $V^0$ is the set of one or more initial points provided by the user.

The search step enables to use various strategies to explore the space of variables. When the search is unsuccessful, i.e., when no better solution is found, a poll step is launched. Every candidate point generated during the poll step must lie within a frame centred around the incumbent solution $x^k$ and which size is parameterized by the \textit{poll size parameter} $\Delta^k\geq\delta^k$.

At the end of an iteration, the mesh and poll size parameters are updated depending on the outcome. If the iteration is unsuccessful both are increased, and conversely, if the iteration is successful both are decreased in such a way that the set of possible directions during the poll gets richer.  In this work, \orthomads~\cite{AbAuDeLe09} is used to deterministically generate $2n$ orthogonal directions at the poll step. A high-level description of \mads is given Algorithm~\ref{algo:MADS}.

%-------------------------%
\begin{algorithm}[htbp]
\caption{The Mesh Adaptive Direct-Search algorithm (\mads).}
\SetAlgoNoLine
\SetKwProg{subalg}{}{}{}

\subalg{\textbf{0. Initialization}}
{\SetAlgoVlined
$\begin{array}{ll}
V^0\subset\mathbb{R}^n &\mbox{: set of starting points}\\
\Delta^0\geq\delta^0>0 &\mbox{: initial mesh and poll size parameters}\\
k\leftarrow 0 &\mbox{: iteration counter}
\end{array}$
}

\subalg{\textbf{1. Search (optional)}}
{\SetAlgoVlined
Evaluate a finite set of points included in the mesh $M^k$.

If the search if successful, go to 3, otherwise go to 2.
}

\subalg{\textbf{2. Poll}}
{\SetAlgoVlined
Evaluate a finite set of points included in the poll frame.
}

\subalg{\textbf{3. Update parameters}}
{\SetAlgoVlined
Update the cache $V^{k+1}$ with the newly sampled points.

Update the mesh and poll size parameters $\delta^{k+1}$ and $\Delta^{k+1}$.

Increase the iteration counter $k\leftarrow k+1$ and go to 1.
}
\label{algo:MADS}
\end{algorithm}
%-------------------------%
%
The algorithm stops either when the poll size parameter falls under a given threshold or when the prescribed budget of function evaluations is spent. Using the Clarke nonsmooth calculus~\cite{Clar83a}, one can prove that under some mild assumptions on the smoothness of the problem, the \mads algorithm globally converges to a solution satisfying local optimality conditions provided that all candidate points lie on the mesh $M^k$. The interested reader may refer to~\cite{AuDe2006} and~\cite[Chapter~8]{AuHa2017}.

In the \mads context, surrogates can be used to find new candidate points during the search step. For instance, the few best solutions of a subproblem that uses only surrogates might be some promising candidate points. Several works have tackled the incorporation of surrogates in \mads like quadratic models~\cite{CoLed2011}, treed Gaussian processes~\cite{GrLeD2011}, LOWESS models~\cite{TaAuKoLed2016}, hybrid models between static surrogates and dynamic models~\cite{AuCM2019}, or ensembles of model~\cite{AuKoLedTa2016}.

%-------------------------------------------%
\section{Quantifying uncertainty with ensembles of models}\label{sec:contribution}
%-------------------------------------------%

Ensembles of models enable to combine several models in the hope of taking advantage of each of them. However, this technique creates an aggregate model that can produce a prediction at any point but not an uncertainty on that prediction, thus prohibiting any Bayesian-like approach. Yet, because several models are useful to describe a single function, it means that their predictions are not identical. Consequently, there should be areas in the search space where the predictions show discrepancies, thus resulting in some form of uncertainty that is not apparent in the aggregate prediction.
The proposed approach is precisely to catch the discrepancies between the models in order to produce a measure of uncertainty.

The idea that the disparity of the models' predictions can be used is tackled in~\cite{GoTuHaShQu2007} where the uncertainty at a given point $x$ is produced with the standard deviation between the predictions defined by
$$\sigma(x) = \left(\frac{\sum_{p=1}^s \big(\tilde f^p(x) - \bar f(x)\big)^2} {s-1}\right)^{\frac{1}{2}}$$
where $\bar f(x) = \sum_{p=1}^s \tilde f^p(x)/s$. This metric quantifies the gaps between the values of the models, however, as it was said earlier, in a BBO context the actual values matter less than the variations of the models. For instance, the two following models possess significantly different values: $\tilde f^1$ and $\tilde f^2=\tilde f^1/10 + 20$. Yet, their variations are the same, i.e., when $\tilde f^1$ increases, $\tilde f^2$ increases too and reciprocally, and therefore they have the same optima. In this case the uncertainty shall be minimum since using either model will yield the same candidate points. Now the two following models differ: $\tilde f^1$ and $\tilde f^3=-\tilde f^1$, but in addition their variations will be opposite so that their optima will certainly be different. In this case, the uncertainty shall be maximum even though the actual standard deviation between $\tilde f^1$ and $\tilde f^3$ might be less than between $\tilde f^1$ and $\tilde f^2$. In light of this, a measure of uncertainty suited to BBO should rather take into account the variations of the models in the form of some local correlation.
In~\cite{Mu2020} and~\cite{To2015}, a correlation coefficient is built between a high-fidelity surrogate $\tilde f^\mathrm{high}$ and a low-fidelity model $\tilde f^\mathrm{low}$
$$   r = \frac{\displaystyle\sum_{j=1}^M \left( \tilde f^\mathrm{high}\big(x^{(j)}\big) - \bar f^\mathrm{high} \right)  \left( \tilde f^\mathrm{low}\big(x^{(j)}\big) - \bar f^\mathrm{low} \right)}{\sqrt{\displaystyle\sum_{j=1}^M \left( \tilde f^\mathrm{high}\big(x^{(j)}\big) - \bar f^\mathrm{high} \right)^2} \sqrt{\displaystyle\sum_{j=1}^M \left( \tilde f^\mathrm{low}\big(x^{(j)}\big) - \bar f^\mathrm{low} \right)^2}}$$
where $\{x^{(j)}\}_{j\in\{1,2,\dots,M\}}$ is a set of $M\geq n$ points sampled locally around an area of interest; and $\bar f^\mathrm{high}$ and $\bar f^\mathrm{low}$ are the average values of $\tilde f^\mathrm{high}$ and $\tilde f^\mathrm{low}$ on this set of points, respectively. For the reasons aforementioned, this quantity is more relevant than the standard deviation in BBO.

%-------------------------------------------------%
\subsection{A new expression for the uncertainty}
%-------------------------------------------------%

In the proposed approach, this idea of correlation between models is exploited to produce an expression of the uncertainty at a given point $x$. Two alternatives are built: a smooth and a nonsmooth uncertainties. In addition each alternative is declined into two versions: an uncertainty dedicated to the objective and another one dedicated to the constraints.

%------------------------------------------------------%
\subsubsection*{Smooth uncertainty for the objective}
%------------------------------------------------------%

The simplex gradient~\cite{CTKelley_1999_b} of a function $f$ at point $x$, denoted by $\nabla_S f(x)$, is the gradient of the linear model of~$f$ at $x$. Computing the simplex gradient around a given point $x$ requires the evaluation of $f$ on a simplex, i.e., a set of $n+1$ affinely independent points, around $x$. The correlation between two models can be reinterpreted geometrically. For two models $\tilde f^p$ and $\tilde f^q$, the cosine between their simplex gradients at a given point $x$ is defined by
$$\cos\left\langle\nabla_S\tilde f^p(x),\ \nabla_S\tilde f^q(x)\right\rangle = \frac{\nabla_S\tilde f^p(x)^{\top}\nabla_S\tilde f^q(x)}{||\nabla_S\tilde f^p(x)||_2\times||\nabla_S\tilde f^q(x)||_2} \mbox{ .}$$
The larger the cosine, the more correlated the models around $x$. With this notion in mind, the uncertainty between two models $\hat\sigma_{p,q}$ can be produced as an inversely proportional function of the cosine
\begin{equation}\label{eq:uncertainty_obj_smooth}
    \hat\sigma_{p,q}(x) := \frac{1}{2}\left(1 - \cos\left\langle\nabla_S\tilde f^p(x),\ \nabla_S\tilde f^q(x)\right\rangle\right)
    \mbox{ .}
\end{equation}
When the models are highly correlated, the cosine is close to 1 so that the uncertainty is close to its minimum 0. When the models are poorly correlated, the cosine is closer to 0 and the uncertainty increases to 0.5. And when the models are anti-correlated, the cosine is close -1 so that the uncertainty reaches its maximum 1. The choice of a simplex is left at the discretion of the user. A small simplex around $x$ will yield a simplex gradient that is a good approximation of the true gradient for smooth functions, but a wider simplex will have a smoothing effect that can be appreciable with nonsmooth or noisy functions. In the current context, the simplex gradient acts as a simple surrogate for the true gradients of the models $\tilde f^p$, $p\in\{1,2,\dots,s\}$ that are not always easy to obtain. However, using the true gradients if available might be equally efficient.
See Appendix~\ref{appendix:pps_simplex} for the practical construction of the simplex used in this work.

The generalization of this expression to more than two models will be described after the other versions of uncertainties are introduced.

%------------------------------------------------------%
\subsubsection*{Nonsmooth uncertainty for the objective}
%------------------------------------------------------%

An alternative for the uncertainty that does not require the computation of simplex gradients is proposed. It requires a positive spanning set of directions $\mathcal{D}$, i.e., a set of at least $n+1$ vectors of $\mathbb{R}^n$ such that any point of $\mathbb{R}^n$ can be written as a positive linear combination thereof~\cite{Davi54b}. The nonsmooth alternative is defined by
\begin{equation}\label{eq:uncertainty_obj_nonsmooth}
    \hat\sigma_{p,q}(x) = \frac{1}{|\mathcal{D}|}\sum_{d\in\mathcal{D}}\mathsf{xor}\Big(\tilde f^p(x+d)<\tilde f^p(x)\ ,\ \tilde f^q(x+d)<\tilde f^q(x)\Big)
\end{equation}
where $\mathsf{xor}(\cdot\ ,\cdot)$ is the exclusive or logical operator. For each direction $d\in\mathcal{D}$, the uncertainty increases if the models predict contradictory trends, i.e., if model $\tilde f^p$ increases from $x$ to $x+d$ while $\tilde f^q$ decreases, or conversely. The term $1/|\mathcal{D}|$ scales the sum between 0 and 1 so that the uncertainty will not be influenced by the number of directions in $\mathcal{D}$. Here again, the choice of a positive spanning set is left at the discretion of the user. Not only can the size of the directions vary as with the simplex, but also the number of directions can increase in order to explore the surroundings of $x$ better.  See Appendix~\ref{appendix:pps_simplex} for the practical construction of the positive spanning set used in this work.

%------------------------------------------------------%
\subsubsection*{Smooth uncertainty for the constraints}
%------------------------------------------------------%

The expressions for the uncertainty proposed in Equations~\eqref{eq:uncertainty_obj_smooth} and~\eqref{eq:uncertainty_obj_nonsmooth} are suited for the objective $f$ since they take into account variations of the models. However, when handling constraints, the key information is the sign of the function rather than whether it increases or not. If two models $\tilde c_j^p$ and $\tilde c_j^q$ of the same constraint $c_j$ are available, the uncertainty shall increase when one model or the other tends towards 0, and increase even more when their signs are opposite, meaning that their predictions on the feasibility are contradictory. Hence the following expression 
\begin{equation}\label{eq:uncertainty_con_smooth}
    \hat\sigma_{p,q}(x) = \mathrm{sigm}\left( -\tilde c_j^p(x)\times\tilde c_j^q(x) \right)
\end{equation}
where sigm$(\cdot)$ is the sigmoid function. It acts as an activation function that increases as the product of $\tilde c_j^p$ and $\tilde c_j^q$ decreases. This uncertainty also ranges from 0 to 1.

%--------------------------------------------------------%
\subsubsection*{Nonsmooth uncertainty for the constraints}
%--------------------------------------------------------%

Here again, a nonsmooth alternative is proposed to the smooth uncertainty for the constraints. It uses the logical operator $\mathsf{xor}$ to indicate whether the two models $\tilde c_j^p$ and $\tilde c_j^q$ predict the same feasibility result at a given point $x$ or not
\begin{equation}\label{eq:uncertainty_con_nonsmooth}
    \hat\sigma_{p,q}(x) = \mathsf{xor}\left(\tilde c_j^p(x)\leq0\ ,\ \tilde c_j^q(x)\leq0\right)
\end{equation}
%

%----------------------------------------------------------------%
\subsubsection*{Generalization to an arbitrary number of models}
%----------------------------------------------------------------%

Expressions~\eqref{eq:uncertainty_obj_smooth},~\eqref{eq:uncertainty_obj_nonsmooth},~\eqref{eq:uncertainty_con_smooth} and~\eqref{eq:uncertainty_con_nonsmooth} consider two models of the objective or a constraint. Ensembles of models usually comprise more than two models though, hence the need for a general expression that can consider an arbitrary number of models. In addition, this general expression must take into account the weights $w^p$, $p\in\{1,2,\dots,s\}$, which reflect the quality of the models. Just as in the prediction defined in Equation~\eqref{eq:aggregated_model}, the good models should have a strong influence in the determination of the uncertainty whereas the poor models should not. The following quantity meets those requirements
\begin{equation}\label{eq:uncertainty_ratio} \left(\sum_{p=1}^{s-1}\sum_{q=p+1}^s w^p w^q\times \hat\sigma_{p,q}(x)\right)\ \Big/\ \sum_{p=1}^{s-1}\sum_{q=p+1}^s w^p w^q
\end{equation}
The ratio in~\eqref{eq:uncertainty_ratio} considers all the possible pairs of models once. For each pair $(p,q)\in\{1,2,\dots,s\}^2$ such that $p\neq q$, the uncertainty $\hat\sigma_{p,q}(x)$ stemming from the models $\tilde f^p$ and $\tilde f^q$ at point $x$ is weighted by the product of the corresponding weights $w^p$ and $w^q$. Consequently, the better the models, the more $\hat\sigma_{p,q}(x)$ will weight up in the total uncertainty at point $x$. Then the sum on all pairs of models is normalized by $1/\sum_{p=1}^{s-1}\sum_{q=p+1}^s w^p w^q$ so that the result does not depend on the number of models~$s$. Since the four versions of $\hat\sigma_{p,q}$ range from 0 to 1, this ratio applies to any case: objective and constraint versions, smooth and nonsmooth alternatives.

At this point the uncertainty takes into account an arbitrary number of models and also the weights as required. But the ratio in~\eqref{eq:uncertainty_ratio} is between 0 and 1 by construction, and therefore it is most likely not at the right scale for the problem at hand. A final step is to multiply by a factor $\alpha>0$ that scales the ratio in a relevant way. For this purpose, $\alpha = 10\times\mathrm{Var}(g(V))$ was chosen, where $g$ is either the objective or a constraint; $g(V)=\left\{ g\big(x^{(1)}\big),g\big(x^{(2)}\big),\dots,g\big(x^{(N_s)}\big) \right\}$ is the set of already sampled values of the function $g$; and 10 is a factor that empirically gave better results. This choice is motivated by the fact that the problem's scale can only be known through the true function's values. The final expression of the uncertainty is 
\begin{equation}\label{eq:uncertainty}
    \hat\sigma(x) = \alpha\frac{w^{\top}\Sigma(x)w}{w^{\top}Tw}
\end{equation}
where the ratio~\eqref{eq:uncertainty_ratio} has been rewritten in a more compact form with $\Sigma(x)\in\mathbb{R}^{s\times s}$ being the upper triangular matrix such that $[\Sigma(x)]_{p,q}=\hat\sigma_{p,q}(x)$ if $p<q$ and 0 otherwise; $T\in\mathbb{R}^{s\times s}$ being the upper triangular matrix such that $[T]_{p,q}=1$ if $p<q$ and 0 otherwise; and $w$ being the vector of weights $[w^1,w^2,\dots,w^s]^{\top}$.

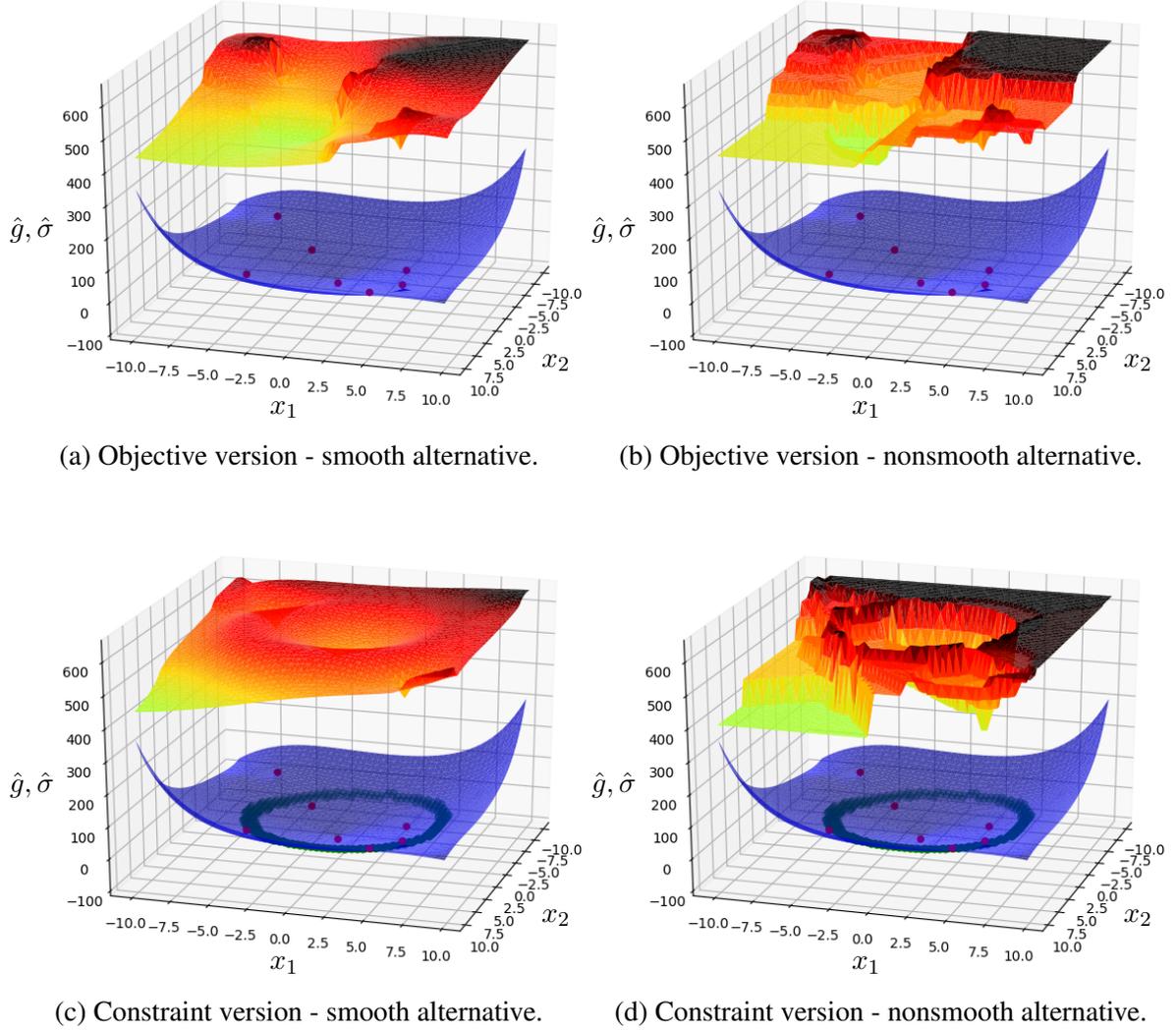
\begin{figure}[!ht]
\setlength{\unitlength}{0.8cm}
\begin{subfigure}{.5\textwidth}
  \centering
  \includegraphics[width=1.1\linewidth]{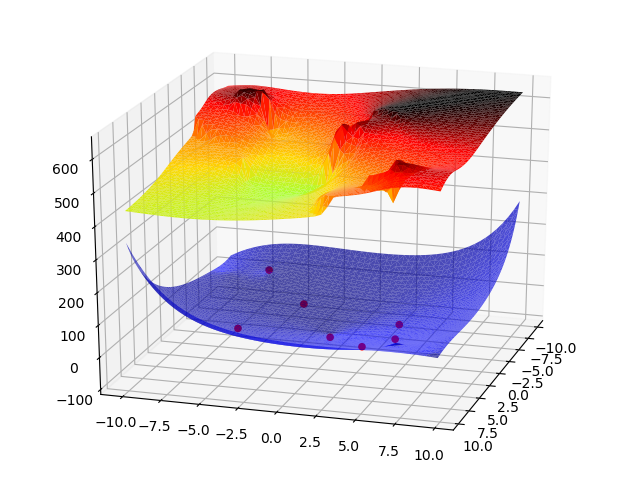}
  \caption{Objective version - smooth alternative.\label{fig:obj_smooth}}
  \begin{picture}(0,0)
  \put(-4.9,4.5){\small$\hat g,\hat\sigma$}
  \put(-0.5,1.5){\small$x_1$}
  \put(4.1,2.3){\small$x_2$}
  \end{picture}
\end{subfigure}
\begin{subfigure}{.5\textwidth}
  \centering
  \includegraphics[width=1.1\linewidth]{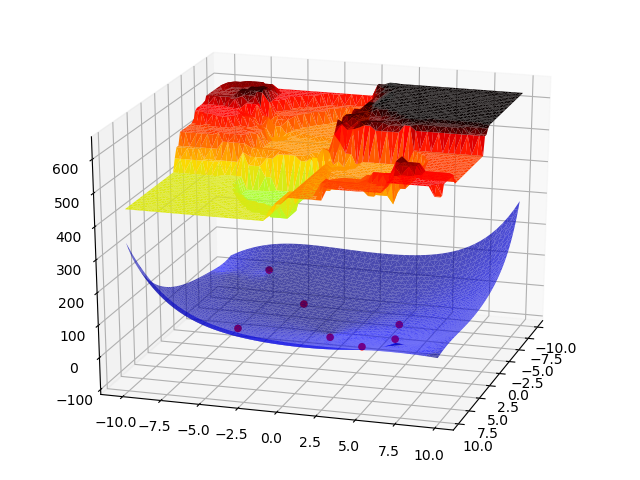}
  \caption{Objective version - nonsmooth alternative.\label{fig:obj_nonsmooth}}
  \begin{picture}(0,0)
  \put(-4.9,4.5){\small$\hat g,\hat\sigma$}
  \put(-0.5,1.5){\small$x_1$}
  \put(4.1,2.3){\small$x_2$}
  \end{picture}
\end{subfigure}
\begin{subfigure}{.5\textwidth}
  \centering
  \includegraphics[width=1.1\linewidth]{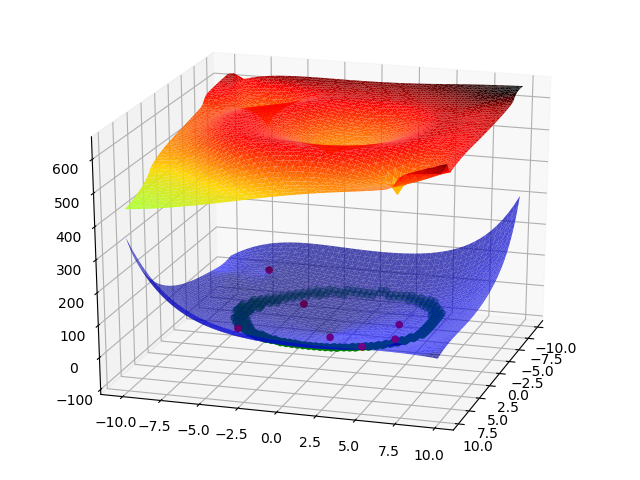}
  \caption{Constraint version - smooth alternative.\label{fig:con_smooth}}
  \begin{picture}(0,0)
  \put(-4.9,4.5){\small$\hat g,\hat\sigma$}
  \put(-0.5,1.5){\small$x_1$}
  \put(4.1,2.3){\small$x_2$}
  \end{picture}
\end{subfigure}
\begin{subfigure}{.5\textwidth}
  \centering
  \includegraphics[width=1.1\linewidth]{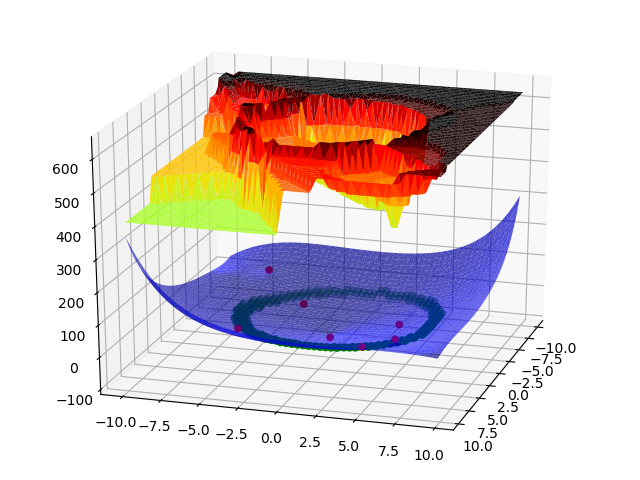}
  \caption{Constraint version - nonsmooth alternative.\label{fig:con_nonsmooth}}
  \begin{picture}(0,0)
  \put(-4.9,4.5){\small$\hat g,\hat\sigma$}
  \put(-0.5,1.5){\small$x_1$}
  \put(4.1,2.3){\small$x_2$}
  \end{picture}
\end{subfigure}

\caption[The four uncertainties on the same sample set]
{The four uncertainties on the same sample set.
Figures~\ref{fig:obj_smooth} and~\ref{fig:obj_nonsmooth} correspond to the the smooth and nonsmooth alternatives of the objective version, respectively (Equations~\eqref{eq:uncertainty_obj_smooth} and~\eqref{eq:uncertainty_obj_nonsmooth}).
Figures~\ref{fig:con_smooth} and~\ref{fig:con_nonsmooth} correspond to the the smooth and nonsmooth alternatives of the constraint version, respectively (Equations~\eqref{eq:uncertainty_con_smooth} and~\eqref{eq:uncertainty_con_nonsmooth}).
}
\label{fig:uncertainty_illustration}

\end{figure}

The different uncertainties are illustrated in
Figure~\ref{fig:uncertainty_illustration}.
Seven points have been sampled in $[-10,10]^2$ from an unknown function $g$ that takes two input variables $x_1$ and $x_2$. On each subfigure, the bottom surface is the aggregate prediction $\hat g$ resulting from eleven different polynomial and RBF models, and the top surface, shifted up for readability, is the uncertainty on that prediction which becomes darker as it increases. The weights of the models $w^p$, $p\in\{1,2,\dots,s\}$, are attributed as described in Section~\ref{subsec:weight_attribution}.
In Figures~\ref{fig:obj_smooth} and~\ref{fig:obj_nonsmooth} the function is interpreted as the objective whereas in
Figures~\ref{fig:con_smooth} and~\ref{fig:con_nonsmooth} the same function is interpreted as a constraint, resulting in significantly different uncertainties. In addition,
in Figures~\ref{fig:con_smooth} and~\ref{fig:con_nonsmooth}, the points close to the assumed border of the constraint, i.e., where the prediction is close to zero, were darkened in order to better understand the uncertainty.

\begin{figure}[!ht]
\setlength{\unitlength}{0.9cm}
    \centering
    \includegraphics[width=0.7\linewidth]{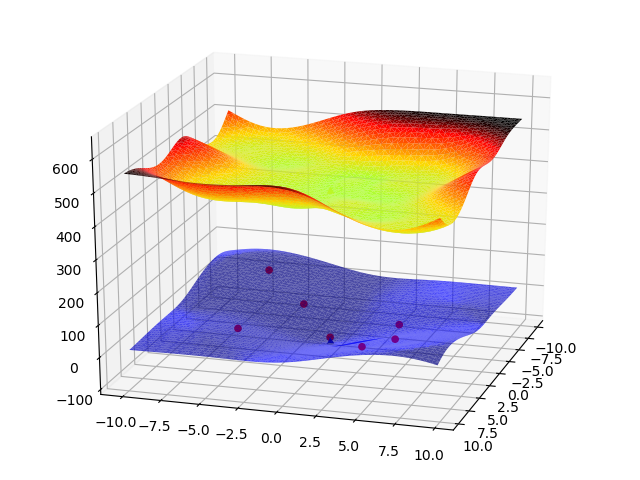}
    \caption{Prediction and uncertainty of Gaussian processes.}
    \label{fig:kriging_uncertainty_illustration}
    
    \begin{picture}(0,0)
    \put(-6,5){$\hat g,\hat\sigma$}
    \put(-1.2,1.6){$x_1$}
    \put(4.2,2.5){$x_2$}
    \end{picture}
\end{figure}

Comparatively, Figure~\ref{fig:kriging_uncertainty_illustration} shows a GP model's prediction and uncertainty fit on the same sample points. It can be noticed that the uncertainty is lower close to the sample points and increases with the distance to them, which is expected with a GP model. The same observation cannot be made in Figure~\ref{fig:uncertainty_illustration}, especially in Figures~\ref{fig:con_smooth} and~\ref{fig:con_nonsmooth} where the uncertainty is higher close to the border of the constraint.
%
%----------------------------------%
%----------------------------------%

%----------------------------------%
%----------------------------------%

%-------------------------------------------------%
\subsection{Error metric and weight attribution}
\label{subsec:weight_attribution}
%-------------------------------------------------%

In this work, the strategy to assign weights is to compute an error metric $\mathcal{E}^p$ for each model $\tilde f^p$ and then attribute a weight $w^p$ that is a function of $\mathcal{E}^p$. The error metric chosen is the order error cross-validation metric (OECV)~\cite{AuKoLedTa2016} mentioned earlier and denoted by $\mathcal{E}_{OECV}$. Broadly speaking, it measures a model capacity to rank points in the same order as the actual objective would do, or, for that matter, its capacity to predict the same feasibility result as an actual constraint would do. When the metric $\mathcal{E}_{OECV}^p$ is computed for every model $\tilde f^p$, $p\in\{1,2,\dots,s\}$, the weights can be attributed. Assigning a weight of 1 to the best model and 0 to the others is the choice made in~\cite{AuKoLedTa2016}. However, in the present work there must be at least two strictly positive weights otherwise the
ratio~\eqref{eq:uncertainty_ratio} is a division by zero and has no meaning. The approach chosen instead is to select the $N_{\mathrm{best}}$ models that have the smallest error metrics and to assign to them a weight proportional to the metric. Formally, if $I\subset \{1,2,\dots,s\}$ is the subset of the selected models indices, then $w^p \propto \mathcal{E}^{\mathrm{tot}}_I-\mathcal{E}^p$ if $p\in I$ and $w^p=0$ otherwise, where $\mathcal{E}^{\mathrm{tot}}_I$ is the total error of the selected models. The weights are then normalized so that $\sum_{p\in I} w^p = 1$. If more than $N_{\mathrm{best}}$ models have an error metric that is equal to the best metric, all of them will be selected and be assigned an equal weight. Preliminary tests showed that $N_{\mathrm{best}}=3$ and $N_{\mathrm{best}}=4$ were appropriate values for the smooth and nonsmooth alternatives, respectively.

\subsection[Incorporation into the MADS algorithm]{Incorporation into the \mads algorithm}

The \mads algorithm offers an important flexibility through its search step. Many works have already included model-based subproblems (SP) into the search step as described earlier.
The proposed implementation falls within this category. At each iteration, a surrogate problem is solved during the search step and the best solution is used as the next candidate point for the true problem. Having said that, there are many ways to design a surrogate SP.

In~\cite{TaLeDKo2014}, eight different formulations for SP are proposed, denoted by SP1 through SP8. They are specifically designed to take advantage of the double information that stochastic models provide, i.e., the prediction and the uncertainty. These formulations tackle general constraints and involve the following statistical measures:
% EI -- expected vs Expected
$$\left\{
\begin{tabular}{rrcl}

Expected improvement~\cite{JoScWe1998}: &
$\mathrm{EI}(x)$ & $=$ & $ \mathbb{E}[\max(f_{min}-f(x),0)]$ \\

Probability of feasibility~\cite{TaLeDKo2014}: &
$\mathrm{P}(x)$ & $=$ & $ \mathbb{P}[c_j(x) \leq 0 ,  j=1,2,\dots,m]$\\

Probability of improvement~\cite{DRJones_2001}: &
$\mathrm{PI}(x)$ & $=$ & $ \mathbb{P}[f_{min}>f(x)]$ \\
\end{tabular}
\right.$$

Then some other measures are derived from these quantities: the \textit{expected feasible improvement} $\mathrm{EFI}(x) = \mathrm{EI}(x)\mathrm{P}(x)$, the \textit{probability of feasible improvement}
$\mathrm{PFI}(x) = \mathrm{PI}(x)\mathrm{P}(x)$ and the \textit{uncertainty on the feasibility} $\mu(x) = 4\mathrm{P}(x)(1-\mathrm{P}(x))$. In the formulations, all these measures are arranged in different ways in order to highlight various properties of the sample points. The eight formulations of SP are given in Appendix~\ref{appendix:SP_formulations}.

With stochastic models, the probability distribution at any point $x$ is known and used to practically compute $\mathrm{EI}(x)$, $\mathrm{P}(x)$ and $\mathrm{PI}(x)$. However, an aggregate model extended with the proposed uncertainty, although inspired by the stochastic modelling philosophy, has no probabilistic foundation. More specifically, there is no cumulative distribution function available
at a given point $\mathbb{P}[g(x)<g_0]$, for all $g_0\in\mathbb{R}$. Consequently, the statistical quantities defined above are not applicable as such. To address this issue, $\mathrm{P}$, $\mathrm{PI}$ and $\mathrm{EI}$ are replaced by substitutes $\widetilde{\mathrm{P}}$, $\widetilde{\mathrm{PI}}$ and $\widetilde{\mathrm{EI}}$ that are inspired from the case when the stochastic model yields at any point $x$ a value that follows a normal distribution $\mathcal{N}\big(\hat y(x), \hat\sigma^2(x)\big)$, which is the case of GPs.

When $\hat c_j(x) \sim \mathcal{N}\big(\hat y_j, \hat\sigma_j^2\big)$, the expression of $\mathrm{P}$ becomes
$$\mathrm{P}(x)=\displaystyle\prod_{j=1}^m\Phi\left( -\frac{\hat y_j}{\hat\sigma_j} \right)$$
where $\Phi$ is the normal cumulative distribution function. Here the product implies that the constraints are assumed to be independent from each other, which might be incorrect but it is the best available approximation in a BBO context. The proposed adaptation is
$$\widetilde{\mathrm{P}}(x)=\displaystyle\prod_{j=1}^m\mathrm{sigm}_{\lambda}\left( -\frac{\hat y_j}{\hat\sigma_j} \right)$$
where $\mathrm{sigm}_{\lambda}$ is the sigmoid function of parameter $\lambda$, i.e., $\mathrm{sigm}_{\lambda}(x)=\mathrm{sigm}(\lambda x)$. Both $\Phi$ and $\mathrm{sigm}_{\lambda}$ tend to 1 when the ratio $-\hat y_j / \hat\sigma_j$ tends to $+\infty$, i.e., either when $\hat y_j$ is highly negative or when the uncertainty is low for whatever negative value of $\hat y_j$, which in both cases means that the constraint $c_j$ is most likely satisfied. They also tend to 0 when the ratio $-\hat y_j / \hat\sigma_j$ tends to $-\infty$, i.e., either when $\hat y_j$ takes high values or when the uncertainty is low for whatever positive value of $\hat y_j$, which in both cases means that the constraint $c_j$ is most likely not satisfied. $\lambda=2$ produces the closest approximation of $\Phi$ but choosing other values enable to control the shape of the function. Preliminary tests showed that $\lambda=3$ and $\lambda=1$ were interesting values for the smooth and nonsmooth alternatives, respectively.

As for PI, when $\hat f(x) \sim \mathcal{N}\big(\hat y, \hat\sigma^2\big)$, the expression becomes
$$\mathrm{PI}(x) = \Phi\left( \frac{f_{min}-\hat y}{\hat\sigma} \right)$$
where $f_{min}$ is the best know value of the objective. The proposed alternative is  
$$\widetilde{\mathrm{PI}}(x) = \mathrm{sigm}_{\lambda}\left( \frac{f_{min}-\hat y}{\hat\sigma} \right)
$$
Here again, $\Phi$ and $\mathrm{sigm}_{\lambda}$ have the same behaviour but the parameter $\lambda$ enables to control the shape of PI. The values chosen for the smooth and nonsmooth alternatives are $\lambda=0.1$ and $\lambda=0.5$, respectively.

Finally, when $\hat f(x) \sim \mathcal{N}\big(\hat y, \hat\sigma^2\big)$, the expression of EI becomes
$$\mathrm{EI}(x) = (f_{min}-\hat y)\Phi\left( \frac{f_{min}-\hat y}{\hat\sigma} \right) + \hat\sigma\phi\left( \frac{f_{min}-\hat y}{\hat\sigma} \right)$$
where $\phi$ is the normal density function. This expression is intimately related to the Gaussian aspect of the model an therefore is not \textit{a priori} suited for non-Gaussian models, let alone models that are not truly stochastic. However, it possesses interesting properties that are independent from the Gaussian nature of the model and that can be seen as essential to the very notion of expected improvement. Firstly, for a fixed $\hat\sigma$, EI decreases in $\hat y$, tends to 0 when $\hat y$ tends to $+\infty$ and is almost proportional to $\hat y$ when $\hat y$ tends to $-\infty$, which is judicious in a minimization context. In addition, when $\hat y$ gets closer to $f_{min}$, EI becomes almost proportional to $\hat\sigma$, meaning that when the prediction does not improve the objective (i.e., $\hat y \simeq f_{min}$) the expected improvement mostly relies on the uncertainty. Then, for a fixed $\hat f$, EI increases in $\hat\sigma$ and is almost proportional to $\hat\sigma$ when $\hat\sigma$ tends to $+\infty$, which is sensible since for a given prediction the higher the uncertainty, the larger the potential improvement. Finally, when $\hat\sigma$ tends to 0, the behavior of EI depends on the values of $f_{min}$ and $\hat y$: if $f_{min}\geq\hat y$, then EI tends to $f_{min} - \hat y$, and if $f_{min}<\hat y$, then EI tends to 0, meaning that when the uncertainty is low, EI mostly relies on the comparison between $f_{min}$ and the prediction $\hat y$. Taking into account these considerations, the proposed adaptation for EI is very close to the actual EI and is defined by
$$\widetilde{\mathrm{EI}}(x) = (f_{min}-\hat y)\ \mathrm{sigm}_{\lambda}\left( \frac{f_{min}-\hat y}{\hat\sigma} \right) + \hat\sigma\gamma\left( \frac{f_{min}-\hat y}{\hat\sigma} \right)$$
where $\gamma(t)=e^{-t^2/2}$. Here $\lambda=1$ was chosen. The functions $\phi$ and $\gamma$ only differ by a factor $1/\sqrt{2\pi}$ and the reason for choosing $\gamma$ instead of $\phi$ is that this factor is no more justified without an actual stochastic model that produces normal distributions. Moreover, preliminary tests showed that the proposed uncertainty seemed on average lower than the uncertainty provided by a kriging model. The terms P, PI and EI were then replaced by $\widetilde{\mathrm{P}}$, $\widetilde{\mathrm{PI}}$ and $\widetilde{\mathrm{EI}}$ in the formulations of SP.

%-------------------------------------------%
\begin{algorithm}[htbp]
\caption{The \mads algorithm with aggregate models.}
\SetAlgoNoLine
\SetKwProg{subalg}{}{}{}

\subalg{\textbf{0. Initialization}}
{\SetAlgoVlined
$\begin{array}{ll}
\mathrm{SP}\in\{\mathrm{SP}_1,\mathrm{SP}_2,\dots,\mathrm{SP}_8\} &\mbox{: surrogate subproblem formulation}\\
\tilde g^1,\tilde g^2,\dots,\tilde g^s &\mbox{: choice of models for the objective and the constraints}\\
V^0\subset\mathbb{R}^n &\mbox{: set of starting points}\\
\Delta^0\geq\delta^0>0 &\mbox{: initial mesh and poll size parameters}\\
k\leftarrow 0 &\mbox{: iteration counter}
\end{array}$
}

\subalg{\textbf{1. Models and weights update}}
{\SetAlgoVlined
Build or update $\tilde f^1,\tilde f^2,\dots,\tilde f^s$ using the values of $f$ in $V^k$\\
Update $w^1,w^2,\dots,w^s$ using the OECV metric for the objective\\
Build or update $\tilde c^1_j,\tilde c^2_j,\dots,\tilde c^s_j$ using the values of $c_j$ in $V^k$, for $j\in\{1,2,\dots,m\}$\\
Update $w^1_j,w^2_j,\dots,w^s_j$ using the OECV metric for the constraints, for $j\in\{1,2,\dots,m\}$\\
}

\subalg{\textbf{2. Search}}
{\SetAlgoVlined
Solve SP to find the best solution $x^k_{SP}$\\
Project $x^k_{SP}$ onto the mesh $M^k$\\
Evaluate the resulting point with the true problem
}

\subalg{\textbf{3. Standard poll}}

\subalg{\textbf{4. Standard parameters update}}

\label{algo:MADS_with_method}
\end{algorithm}
%-------------------------------------------%

Algorithm~\ref{algo:MADS_with_method} summarizes the incorporation of extended aggregate models in \mads. First, one formulation must be chosen among \{SP1, SP2, $\dots$, SP8\}. At iteration $k$, the best solution found for SP, denoted by $x^k_{SP}$, is projected onto the mesh $M^k$ and is used as the candidate point of the search step. The freshly evaluated points are then added to the cache $V^k$ so that the models and the weights will be adjusted accordingly before iteration $k+1$ begins. The resulting algorithm benefits from the convergence results of \mads since all the candidate points lie on the mesh $M^k$.

%-------------------------------------------%
\section{Computational results}
\label{sec:results}
%-------------------------------------------%

The proposed approach has been tested on seven analytical problems; two multi-disciplinary optimization (MDO) applications: the aircraft range problem and the simplified wing problem; and two simulation problems: \solar{1} and \styrene. Version~4 of the \nomad software~\cite{AuCo04a,nomad4paper} was used on a PC Intel(R) Core(TM) i7-8700~CPU~@~3.20GHz on Linux. The aggregate models used the default selection of eighteen models comprised of polynomial response surfaces of various degrees, kernel smoothing, modified radial basis functions as in~\cite{AuKoLedTa2016}, and closest neighbours. The competing quadratic and kriging models were also readily available in \nomad. Due to the long running times required by the kriging models, the latter were only tested on the analytical problems and the aircraft range problem. Every instance of \mads in this work uses the last direction of success at the poll step~\cite{AuDe2006}. Two other BBO solvers have been included in this study: \dfn~\cite{FaLiLuRi2014} and \shebo~\cite{Mueller2019}. The former innately handles general constraints whereas the latter is designed for problems with hidden constraints. Consequently, the problems were adapted in \shebo so that a violated general constraint will be interpreted as a hidden constraint.

The interpretation of the results mostly relies on data profiles~\cite{MoWi2009} which enable to compare multiple solvers on a given set of problems. Broadly speaking, the data profile of a solver indicates the proportion of problems solved to a given tolerance within a prescribed number of evaluations. Since \shebo does not take a single starting point as an input, it is not fit for comparison to the other algorithms through data profiles. Section~\ref{subsec:results_shebo} provides tabular comparisons with \shebo.

Unless otherwise specified, due to the randomness contained in \mads, every version thereof was run four times on each problem with a different seed for the random generator each time. Similarly, \dfn enables to choose between the Halton and Sobol sequences so the two were tested on each problem and taken into account in the data profiles.

%-----------------------------------%
\subsection{Analytical problems}
\label{subsec:results_analytical}
%-----------------------------------%

The seven analytical problems are listed in
Table~\ref{tab:analytical_problem} with the number of variables $n$ and constraints $m$, whether the variables are bounded or not, and the number of starting points used. By taking into account the additional starting points for problems HS83, HS114 and MAD6, the total number of problems is fifteen. The evaluation budget is $1200(n+1)$.

%-----------------------------------%
\begin{table}[htbp]
\begin{center}
\renewcommand{\tabcolsep}{3pt}
\begin{tabular}{|lllllcc|}
\hline
\# & Name & Source & $n$ & $m$ & Bounds & \# starting points \\ % & Smth & $f^*$ \\
\hline
\hline
 1 &  G2		    &\cite{AuDeLe07} & $10$ & $2$  & yes & 1 \\
 2 &  HS19		    &\cite{HoSc1981} & $2$  & $2$  & yes & 1 \\
 3 &  HS83 		    &\cite{HoSc1981} & $5$  & $6$  & yes & 4 \\
 4 &  HS114		    &\cite{LuVl00}   & $9$  & $4$  & yes & 3 \\
 5 &  MAD6	       	&\cite{LuVl00}   & $5$  & $7$  & no  & 4 \\
 6 &  PENTAGON\ \ \ &\cite{LuVl00}   & $6$  & $15$ & no  & 1 \\
 7 &  SNAKE 	    &\cite{AuDe09a}  & $2$  & $2$  & no  & 1 \\
\hline
\end{tabular}
\end{center}
\caption{Description of the seven analytical problems.}
\label{tab-pbs}
\label{tab:analytical_problem}
\end{table}
%-----------------------------------%

As in~\cite{TaLeDKo2014}, the eight SP formulations were compared, and the following values were tested for the parameter $\lambda$ when applicable: $\{0,0.01,0.1,1\}$, thus resulting in twenty-three distinct formulations. When this parameter is involved in a formulation, it is denoted as a subscript, e.g., SP$2_{0.1}$. The purpose here is not to exhaustively compare the formulations with each other but rather to identify the bests formulations and compare their performances to the existing versions of \nomad and to \dfn.

The formulations with extended aggregate models were compared to \nomad without any search step, referred to as ``\textsf{no search}", \nomad with a search step involving the minimization of quadratic models, referred to as ``\textsf{quad search}", and \dfn. It turns out that all formulations perform better than \textsf{no search}, confirming that the approach is valid and does not ``waste" evaluations. However, \textsf{quad search} is most of the times as good as, and sometimes better than, the proposed approach, confirming the effectiveness of quadratic models on analytical functions. \dfn presents heterogeneous performances. It found good solutions for three problems but performed poorly on the others, hence the low overall performance.

\begin{figure}[htbp]
\begin{subfigure}{.5\textwidth}
  \centering
  \includegraphics[width=\linewidth]{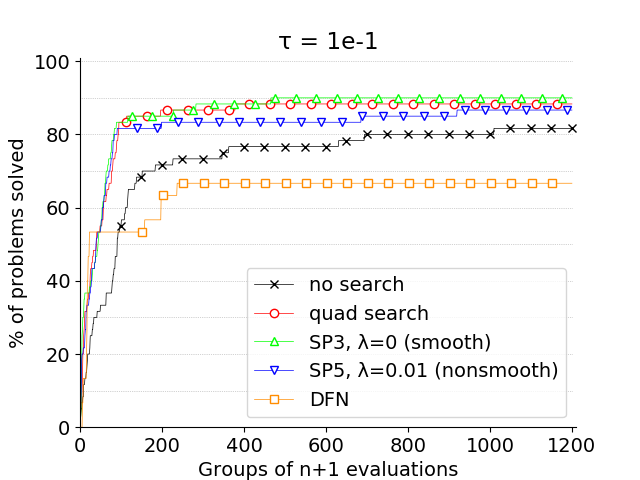}
\end{subfigure}
\begin{subfigure}{.5\textwidth}
  \centering
  \includegraphics[width=\linewidth]{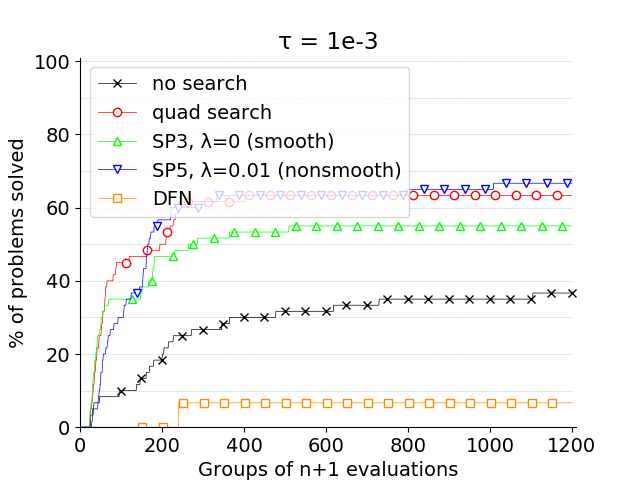}
\end{subfigure}
\begin{subfigure}{.5\textwidth}
  \centering
  \includegraphics[width=\linewidth]{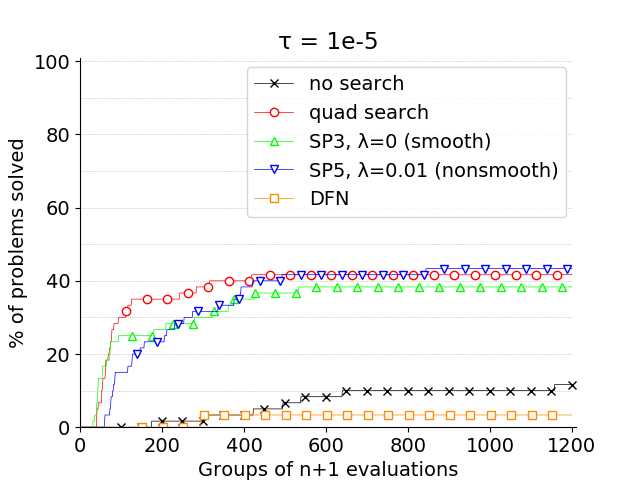}
\end{subfigure}
\begin{subfigure}{.5\textwidth}
  \centering
  \includegraphics[width=\linewidth]{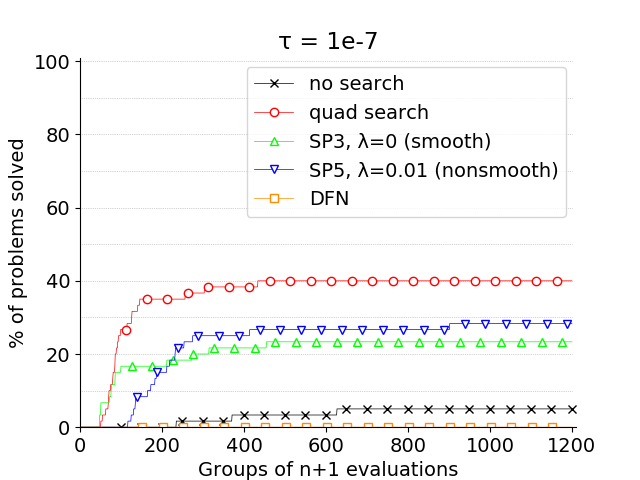}
\end{subfigure}
\caption{Data profiles. \textsf{no search} vs. \textsf{quad search} vs. SP$3_{0}$ with smooth uncertainty vs. SP$5_{0.01}$ with nonsmooth uncertainty vs. \dfn on analytical problems.}
\label{fig:no_search_vs_quadsearch_vs_SP3smooth0_vs_SP5nonsmooth0.01_vs_DFN}
\end{figure}

The best formulation found for the smooth alternative in terms of the proportion of problems solved within the budget is SP$3_{0}$ defined by 
\begin{equation*}
\tag{SP3-EI$\sigma$}
\begin{aligned}
    \min_{x\in\mathcal{X}}\ -&\mathrm{EI}(x) \\
    \mathrm{s.t.}\quad \ &\hat c_j(x) \leq 0,\ \ j=1,2,\dots,m
\end{aligned}
\end{equation*}
As for the nonsmooth alternative, the best formulation is SP$5_{0.01}$ which consists in maximizing $\mathrm{EFI}(x)+0.01\hat\sigma_f(x)$. Figure~\ref{fig:no_search_vs_quadsearch_vs_SP3smooth0_vs_SP5nonsmooth0.01_vs_DFN} shows the data profiles of the five following algorithms: \textsf{no search}, \textsf{quad search}, SP$3_{0}$ with smooth uncertainty, SP$5_{0.01}$ with nonsmooth uncertainty, and \dfn at variable tolerance: $\tau=10^{-1}$, $\tau=10^{-3}$, $\tau=10^{-5}$ and $\tau=10^{-7}$.  

%---------------------%
%---------------------%

The profiles show that the performances of SP3 and SP5 are close to that of \textsf{quad search} for high tolerance. However, \textsf{quad search} becomes significantly better for low tolerance ($\tau=10^{-7}$). These first results show that the extended aggregate models combined with the formulations manage to find good solutions as efficiently as \textsf{quad search}. However, quadratic models do so slightly faster in terms of the number of evaluations, and more importantly, they are especially accurate on analytical problems, thus resulting in superior performances at low tolerance.

% Algorithm \textsf{no search} did not find feasible solution for problem HS19 for seeds 1,2,3. Algo SP8 nonsmooth did not find feasible solution for problem HS19 for seeds 1,3.

\begin{figure}[htbp]
\begin{subfigure}{.5\textwidth}
  \centering
  \includegraphics[width=\linewidth]{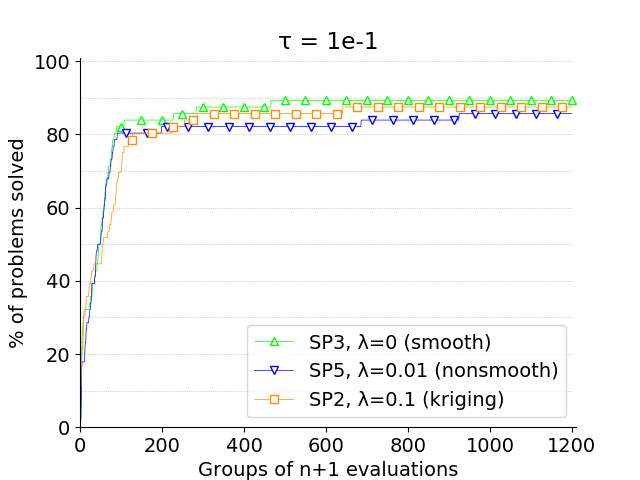}
\end{subfigure}
\begin{subfigure}{.5\textwidth}
  \centering
  \includegraphics[width=\linewidth]{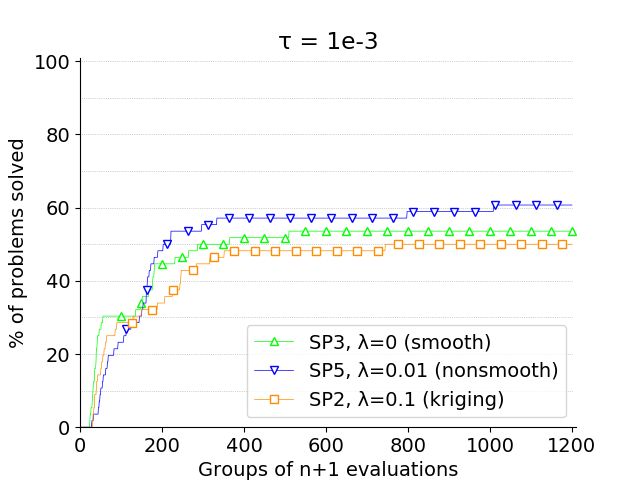}
\end{subfigure}
\begin{subfigure}{.5\textwidth}
  \centering
  \includegraphics[width=\linewidth]{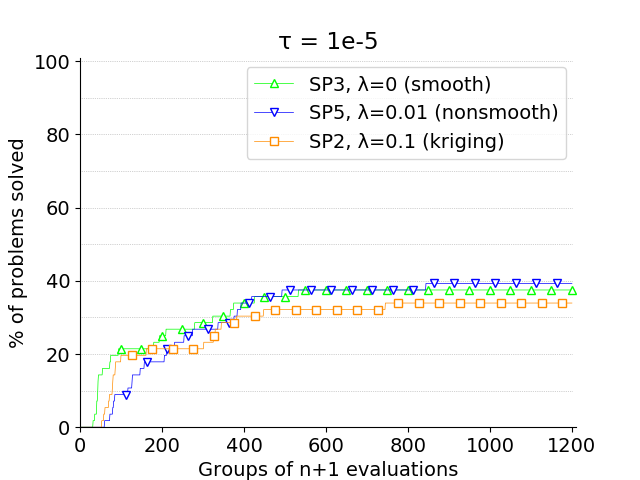}
\end{subfigure}
\begin{subfigure}{.5\textwidth}
  \centering
  \includegraphics[width=\linewidth]{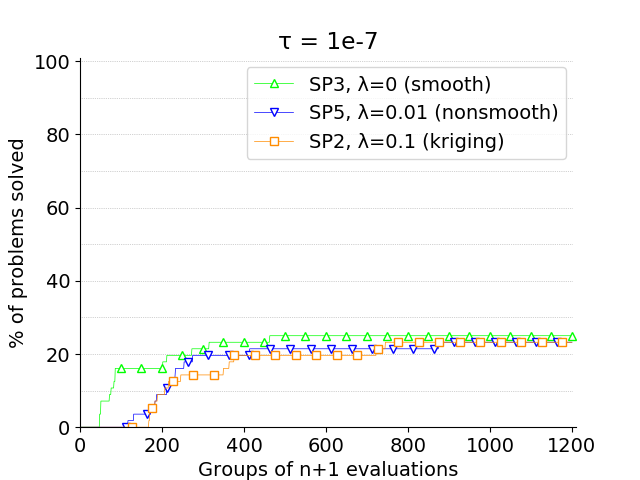}
\end{subfigure}
\caption{Data profiles. SP$3_{0}$ with smooth uncertainty vs. SP$5_{0.01}$ with nonsmooth uncertainty vs. SP$2_{0.1}$ with kriging models on analytical problems.}
\label{fig:SP3smooth0_vs_SP5nonsmooth0.01_vs_SP2kriging0.1}
\end{figure}

Since extended aggregate models are meant to mimic and therefore supersede actual stochastic models, the comparison to the available kriging models in \nomad is most appropriate. In~\cite{TaLeDKo2014}, the authors recommend SP1 and SP2 with large values of $\lambda$ when using stochastic models. For that reason, these two formulations have been tested with $\lambda=0.1$ and $\lambda=1$. On this set of problems, SP$2_{0.1}$ turns out to be the best formulation. The latter was therefore tested against SP$3_{0}$ with smooth uncertainty and SP$5_{0.01}$ with nonsmooth uncertainty, that is the best formulations seen above. The resulting data profiles in Figure~\ref{fig:SP3smooth0_vs_SP5nonsmooth0.01_vs_SP2kriging0.1} show that on the present set of problems extended aggregate models are as good as, or better than, the kriging alternative depending on the tolerance. In addition, due to the inversion of a covariance matrix that grows with the size of the sample set, the kriging models typically take minutes to tens of minutes to solve one problem, which is prohibitive when optimizing cheap functions with large budgets of evaluations. In comparison, the proposed approach and \textsf{quad search} typically take minutes and \textsf{no search} and \dfn take seconds. As a result, replacing classic kriging models by extended aggregate models does not harm the performances in terms of the number of evaluations on the set of analytical problems, while improving the real optimization time. Overall, based on the present results, \textsf{quad search} must be favoured on cheap analytical problems.

%-----------------------------------%
\subsection{The aircraft range MDO problem}
%-----------------------------------%

The aircraft range problem~\cite{Kodiyalam2001} is a multi-disciplinary optimization problem (MDO), meaning that it is combined of several interconnected disciplines so that the input of one discipline is the output of the others, and several cycles between them are necessary to stabilize the result. The aircraft range problem aims at maximizing the range of a supersonic business jet by considering aerodynamics, structure, and propulsion, under constraints of engine, performance and structure. The problem has $n=10$ variables and $m=10$ constraints. It is nonsmooth and has several local optima.

In order to test the algorithms truthfully, ten starting points were sampled with Latin hypercube sampling~\cite{McCoBe79a}. Here again, the 23 formulations were run with the two uncertainty alternatives in order to identify the best combinations and analyze their performances. The evaluation budget is $1000(n+1)$. The best formulation with the smooth alternative is SP$3_{0.1}$ defined by
\begin{equation*}
\tag{SP3-EI$\sigma$}
\begin{aligned}
    \min_{x\in\mathcal{X}}\ -\mathrm{EI}(x)-&0.1\hat\sigma_f(x) \\
    \mathrm{s.t.}\quad \ \hat c_j(x)-&0.1\hat\sigma_j(x)\leq0,\ \ j=1,2,\dots,m
\end{aligned}
\end{equation*}
and the best formulation identified with the nonsmooth alternative is SP8 which consists in maximizing $\mathrm{PFI}(x)$.

\begin{figure}[!ht]
\begin{subfigure}{.5\textwidth}
  \centering
  \includegraphics[width=\linewidth]{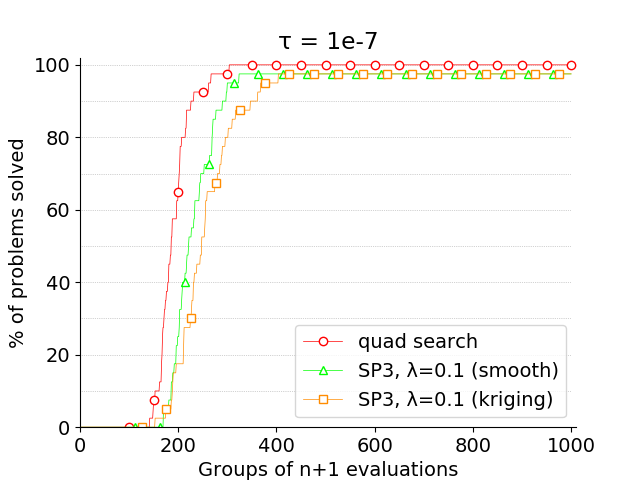}
  \caption{}
  \label{fig:quadsearch_vs_SP3smooth_vs_SP3kriging_dp7}
\end{subfigure}
\begin{subfigure}{.5\textwidth}
  \centering
  \includegraphics[width=\linewidth]{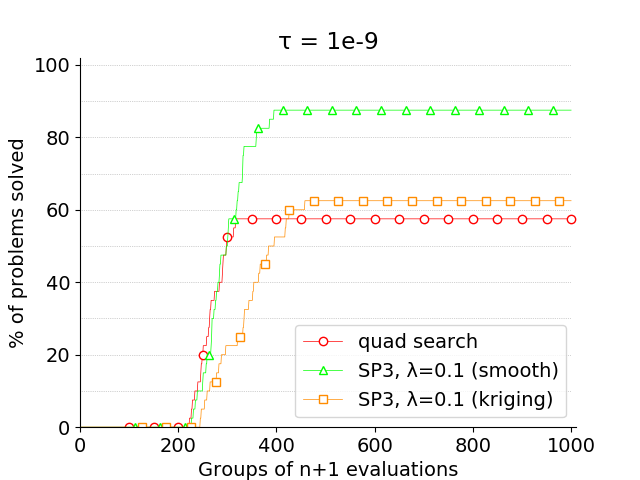}
  \caption{}
  \label{fig:quadsearch_vs_SP3smooth_vs_SP3kriging_dp9}
\end{subfigure}
\begin{subfigure}{.5\textwidth}
  \centering
  \includegraphics[width=\linewidth]{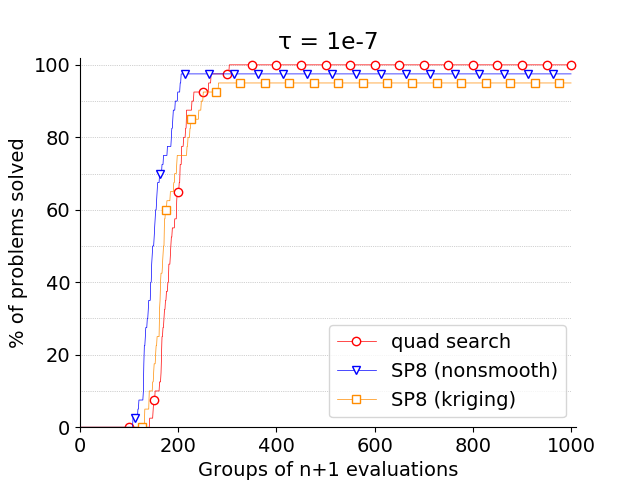}
  \caption{}
  \label{fig:quadsearch_vs_SP8nonsmooth_vs_SP8kriging_dp7}
\end{subfigure}
\begin{subfigure}{.5\textwidth}
  \centering
  \includegraphics[width=\linewidth]{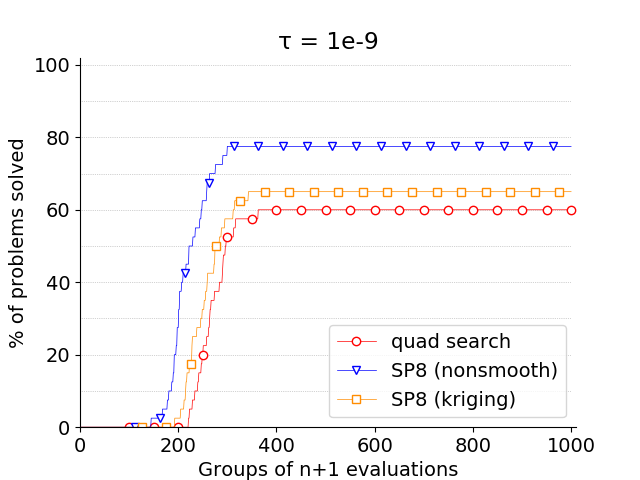}
  \caption{}
  \label{fig:quadsearch_vs_SP8nonsmooth_vs_SP8kriging_dp9}
\end{subfigure}
\caption{Data profiles. Figures~\ref{fig:quadsearch_vs_SP3smooth_vs_SP3kriging_dp7}
and~\ref{fig:quadsearch_vs_SP3smooth_vs_SP3kriging_dp9} show \textsf{quad search} vs. SP$3_{0.1}$ with smooth uncertainty vs. SP$3_{0.1}$ with kriging models.
Figures~\ref{fig:quadsearch_vs_SP8nonsmooth_vs_SP8kriging_dp7}
and~\ref{fig:quadsearch_vs_SP8nonsmooth_vs_SP8kriging_dp9} show \textsf{quad search} vs. SP8 with nonsmooth uncertainty vs. SP8 with kriging models on the aircraft range problem.}
\end{figure}

This time, each of the aforementioned formulations is compared individually not only to \textsf{quad search}, but also to the same formulation with kriging models instead. This choice is motivated by the fact that the aircraft range problem is a real-world simulation-based problem that takes significantly more time to compute than the previous analytical problems, and therefore computing expensive kriging models might be worth the trade-off. In this section, \textsf{no search} is not shown in the interest of readability, and \dfn neither because it performed poorly and its data profiles were flat. Figures~\ref{fig:quadsearch_vs_SP3smooth_vs_SP3kriging_dp7}
and~\ref{fig:quadsearch_vs_SP3smooth_vs_SP3kriging_dp9} show the data profiles of SP$3_{0.1}$ with smooth uncertainty with tolerances $\tau=10^{-7}$ and $\tau=10^{-9}$, respectively; and Figures~\ref{fig:quadsearch_vs_SP8nonsmooth_vs_SP8kriging_dp7}
and~\ref{fig:quadsearch_vs_SP8nonsmooth_vs_SP8kriging_dp9} show the data profiles of SP8 with nonsmooth uncertainty with tolerances $\tau=10^{-7}$ and $\tau=10^{-9}$, respectively. All the algorithms presented equivalent performances for high tolerances ($\tau=10^{-1}$ and $\tau=10^{-3}$) and consequently the data profiles did not reveal significant difference between the solvers with tolerances under $\tau=10^{-5}$, meaning that all the formulations manage to reach a good solution.

For every above-mentioned formulation, the performance is relatively comparable to that of \textsf{quad search} for tolerance $\tau=10^{-7}$. However, for $\tau=10^{-9}$, the proposed extended aggregate models coupled with the right formulations turn out to be significantly better than both \textsf{quad search} and their kriging counterpart, i.e., the same formulations with kriging models instead. It can be noticed that the smooth alternative eventually solves more problem, but the nonsmooth one solves problems faster. In Figure~\ref{fig:quadsearch_vs_SP3smooth_vs_SP3kriging_dp9}, SP$3_{0.1}$ with smooth uncertainty solves 87.5\% of the problems at tolerance $\tau=10^{-9}$ with the allocated budget, while in Figure~\ref{fig:quadsearch_vs_SP8nonsmooth_vs_SP8kriging_dp9} SP8 with nonsmooth uncertainty only solves 77.5\% of the problems at the same tolerance. However, the latter takes only $300(n+1)$ evaluations to do so, while the former has only solved 50\% of the problems after the same number of evaluations. This trend has been observed on this problem with all the other formulations not presented in this paper. The higher achievements of the smooth uncertainty over the long run may be attributed to its intrinsically rich range of values, while the relative rapidity of the nonsmooth uncertainty may be the result of a more aggressive behaviour that helps find a good solution faster.

\begin{figure}[!ht]
\begin{subfigure}{.5\textwidth}
  \centering
  \includegraphics[width=\linewidth]{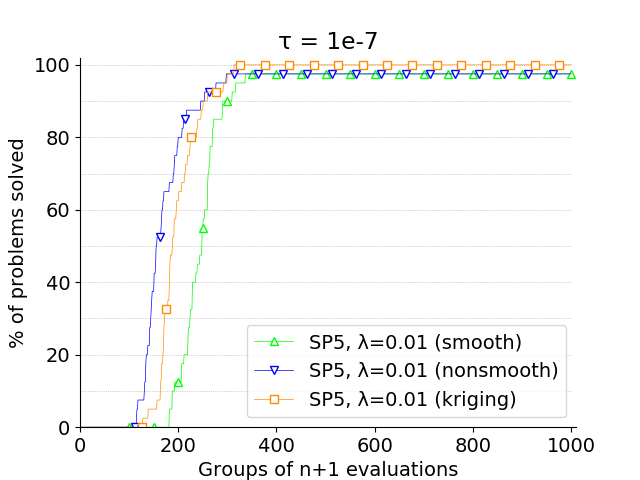}
  \caption{}
  \label{fig:SP5smooth_vs_SP5nonsmooth_vs_SP5kriging_dp7}
\end{subfigure}
\begin{subfigure}{.5\textwidth}
  \centering
  \includegraphics[width=\linewidth]{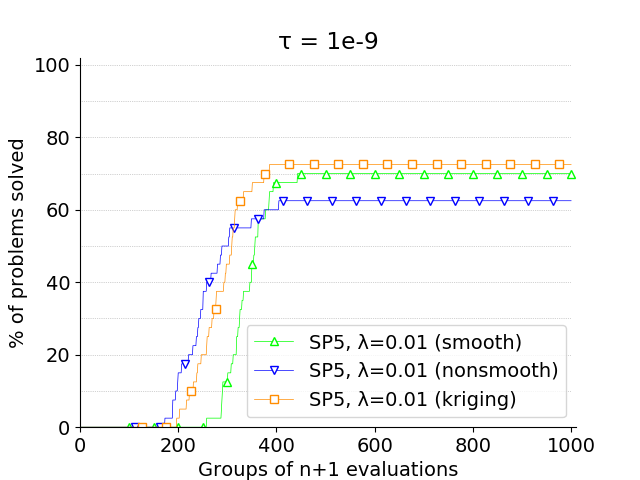}
  \caption{}
  \label{fig:SP5smooth_vs_SP5nonsmooth_vs_SP5kriging_dp9}
\end{subfigure}
\begin{subfigure}{.5\textwidth}
  \centering
  \includegraphics[width=\linewidth]{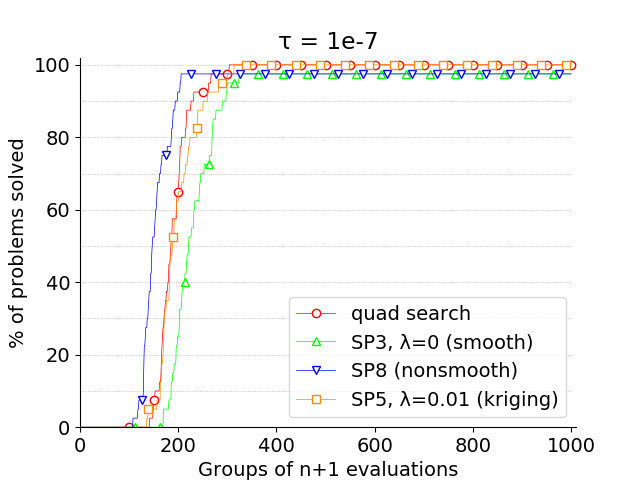}
  \caption{}
  \label{fig:quadsearch_vs_SP3smooth_vs_SP8nonsmooth_vs_SP5kriging_dp7}
\end{subfigure}
\begin{subfigure}{.5\textwidth}
  \centering
  \includegraphics[width=\linewidth]{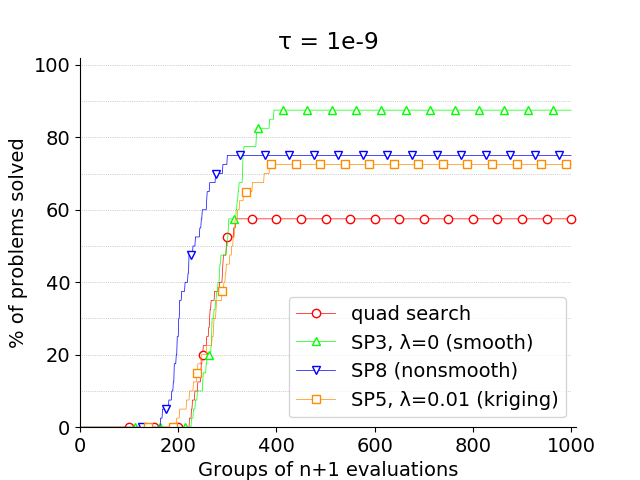}
  \caption{}
  \label{fig:quadsearch_vs_SP3smooth_vs_SP8nonsmooth_vs_SP5kriging_dp9}
\end{subfigure}
\caption{Data profiles. Figures~\ref{fig:SP5smooth_vs_SP5nonsmooth_vs_SP5kriging_dp7} 
and~\ref{fig:SP5smooth_vs_SP5nonsmooth_vs_SP5kriging_dp9} show SP$5_{0.01}$ with smooth uncertainty vs. SP$5_{0.01}$ with nonsmooth uncertainty vs. SP$5_{0.01}$ with kriging models. Figures~\ref{fig:quadsearch_vs_SP3smooth_vs_SP8nonsmooth_vs_SP5kriging_dp7}
and~\ref{fig:quadsearch_vs_SP3smooth_vs_SP8nonsmooth_vs_SP5kriging_dp9} show \textsf{quad search} vs. SP$3_{0.1}$ with smooth uncertainty vs. SP8 with nonsmooth uncertainty vs. SP$5_{0.01}$ with kriging models on the aircraft range problem.}
\end{figure}

It could be argued that comparing kriging models within the best formulations found for extended aggregate models, i.e., SP$3_{0.1}$ and SP8, is not fair since the former might perform poorly on these formulations but yield better results on others. In~\cite{TaLeDKo2014}, the authors use stochastic models and recommend SP5, SP6 and SP7 with small values of $\lambda$ for expensive simulation-based problems as is the case with the aircraft range problem. Accordingly, those three formulations were tested with kriging models and $\lambda=0.01$. On this particular problem, SP5 turns out to be the best formulation. Consequently, extended aggregate models were compared to kriging models on SP5. Figures~\ref{fig:SP5smooth_vs_SP5nonsmooth_vs_SP5kriging_dp7}
and~\ref{fig:SP5smooth_vs_SP5nonsmooth_vs_SP5kriging_dp9} show the data profiles of SP$5_{0.01}$ with smooth and nonsmooth alternatives, and kriging models. The latter perform indeed better than the extended aggregate models on SP$5_{0.01}$ at tolerances $\tau=10^{-7}$ and $\tau=10^{-9}$, meaning that the best formulations with real stochastic models are not necessarily the same than the best ones with the proposed extended aggregate models.

As a result, the fair comparison is not between kriging models and extended aggregate models \textit{within the same formulation}, but rather between each type of models coupled with its best formulation, that is SP$3_{0.1}$ with the smooth alternative, SP8 with the nonsmooth alternative, and SP5 with kriging models.
Figures~\ref{fig:quadsearch_vs_SP3smooth_vs_SP8nonsmooth_vs_SP5kriging_dp7}
and~\ref{fig:quadsearch_vs_SP3smooth_vs_SP8nonsmooth_vs_SP5kriging_dp9} show the data profiles of the above combinations and in addition \textsf{quad search} as a reference. At tolerance $\tau=10^{-7}$, the performances are alike but at tolerance $\tau=10^{-9}$, the extended aggregate models coupled with the suitable formulations solve more problems, faster than kriging models with their own appropriate formulation, and faster than \textsf{quad search}. To tolerance $\tau=10^{-9}$, the smooth and nonsmooth alternatives solve 87.5\% and 75\% of the problems, respectively, while kriging models and \textsf{quad search} solve 72.5\% and 57.5\% of the problems, respectively.

\begin{figure}[!ht]
\begin{subfigure}{.5\textwidth}
  \centering
  \includegraphics[width=\linewidth]{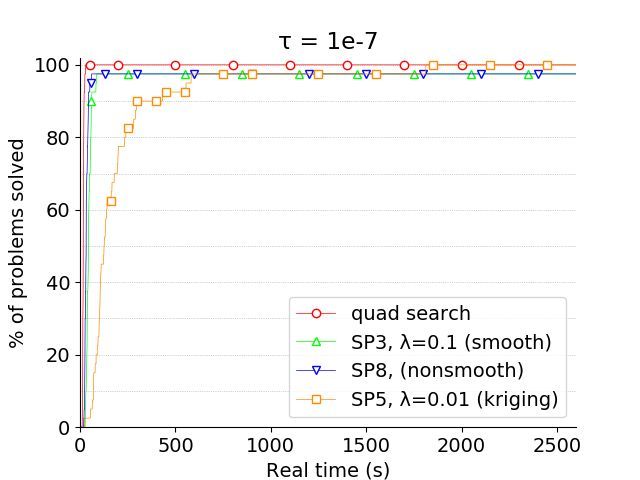}
\end{subfigure}
\begin{subfigure}{.5\textwidth}
  \centering
  \includegraphics[width=\linewidth]{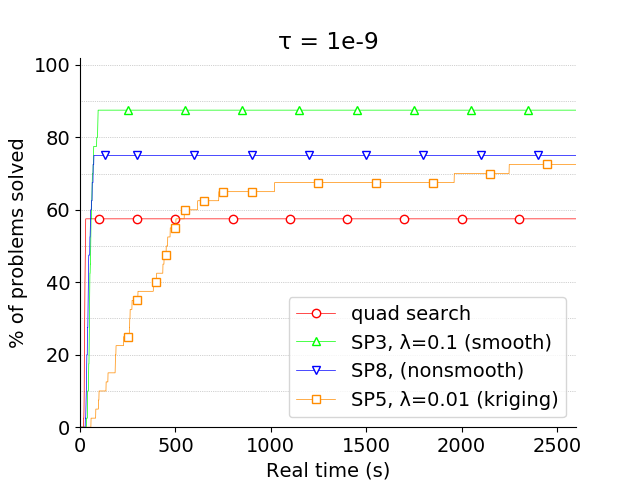}
\end{subfigure}
\caption{Time data profiles. \textsf{quad search} vs. SP$3_{0.1}$ with smooth uncertainty vs. SP8 with nonsmooth uncertainty vs. SP$5_{0.01}$ with kriging models on the aircraft range problem.}
\label{fig:quadsearch_vs_SP3smooth_vs_SP8nonsmooth_vs_SP5kriging_time}
\end{figure}

Since the different approaches require some non negligible amount of internal computation, it is appropriate to compare them not only in terms of the number of evaluations, but also in terms of the total real optimization time. It was highlighted in
Section~\ref{subsec:results_analytical} that kriging models were prohibitively long to train when optimizing cheap analytical problems. The same question is addressed more thoroughly with the aircraft range problem.
Figure~\ref{fig:quadsearch_vs_SP3smooth_vs_SP8nonsmooth_vs_SP5kriging_time} shows the \textit{time data profiles} of formulation SP$3_{0.1}$ with the smooth alternative, formulation SP8 with the nonsmooth alternative, and formulation SP$5_{0.01}$ with kriging models. In a time data profile, the proportion of problems solved is not a function of the number of evaluations anymore but a function of the real computation time instead. The profiles suggest that even with a more expensive, real-world problem, kriging models are especially long to train. At tolerance $\tau=10^{-9}$, SP$3_{0.1}$ with the smooth alternative solves 87.5\% of the problems within 95 seconds and SP8 with the nonsmooth alternative solves 75\% of the problems within 73 seconds. After the same amount of time, SP5 with kriging models has solved less than 10\% of the problems, and requires 728 seconds to solve 65\% of the problems. Some instances even require more than 1500 seconds. Based on the present results, the proposed extended aggregate models constitute a cheaper and more efficient alternative to kriging models.

%-----------------------------------%
\subsection{The simplified wing problem}
%-----------------------------------%

The simplified wing problem~\cite{TriDuTre04a} is also an MDO problem. It aims at minimizing the drag of a wing by taking into account aerodynamics and structure. This problem is smooth but has several local optima. It has $n=7$ bounded variables and $m=3$ constraints.

\begin{figure}[!ht]
\begin{subfigure}{.5\textwidth}
  \centering
  \includegraphics[width=\linewidth]{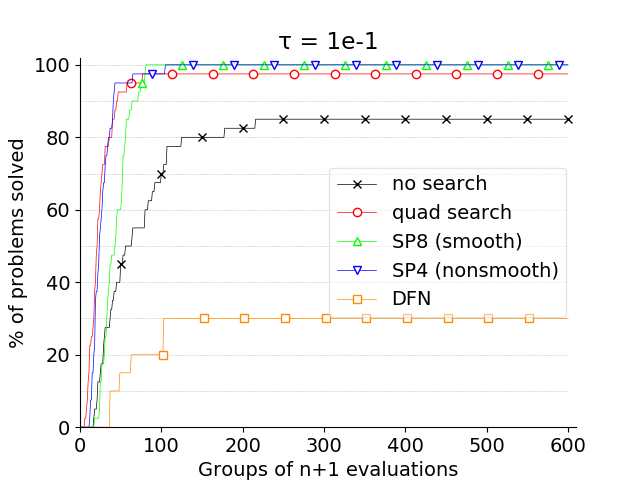}
\end{subfigure}
\begin{subfigure}{.5\textwidth}
  \centering
  \includegraphics[width=\linewidth]{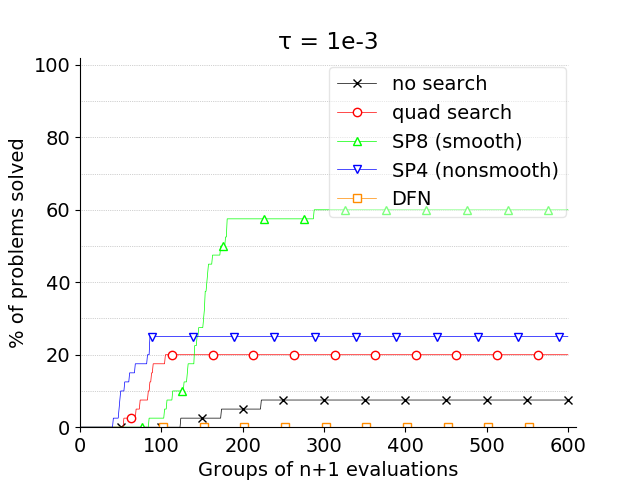}
\end{subfigure}
\caption{Data profiles. \textsf{no search} vs. \textsf{quad search} vs. SP8 with smooth uncertainty vs. SP4 with nonsmooth uncertainty vs. DFN on the simplified wing problem.}
\label{fig:nosearch_vs_quadsearch_vs_SP8smooth_vs_SP4nonsmooth_vs_DFN}
\end{figure}

Ten starting point have been randomly sampled in the bounded space of variables with Latin hypercube sampling. The evaluation budget is $600(n+1)$. Unlike the previous problems and for the rest of the study, the kriging models have not been tested because of their heavy computational cost, and not all the formulations have been tested but rather a subset of the most promising ones found on the aircraft range problem: SP$1_{1}$, SP$2_{0}$, SP$3_{0.1}$, SP$7_{1}$ and SP8 for the smooth alternative; and SP$2_{0.1}$, SP$3_{0.01}$, SP4, SP$6_{0.1}$ and SP8 for the nonsmooth alternative. Among the above formulations, the best ones found for this problem with the smooth and nonsmooth alternatives are SP8 and SP4, respectively. SP4 consists in maximizing $\mathrm{EFI}(x)$. Figure~\ref{fig:nosearch_vs_quadsearch_vs_SP8smooth_vs_SP4nonsmooth_vs_DFN} shows the data profiles of \textsf{no search}, \textsf{quad search}, SP8 with smooth uncertainty, SP4 with nonsmooth uncertainty, and DFN at tolerances $\tau=10^{-1}$ and $\tau=10^{-3}$. Lower tolerances resulted in data profiles that were too flat, especially because of the high sensitivity to the seed of \nomad on this problem.

On this MDO problem, the extended aggregated models perform better than \textsf{no search}, \textsf{quad search} and \dfn. SP4 with the smooth uncertainty is especially good at tolerance $\tau=10^{-3}$, solving 60\% of problems while \textsf{quad search} only solves 20\%. As with the aircraft range problem, the nonsmooth alternative solves problems faster than the smooth one. However, on the simplified wing problem the smooth alternative is by far the most efficient eventually.

%-----------------------------------%
\subsection[The solar{1} problem]{The \solar{1} problem}
%-----------------------------------%

The \solar{1} problem is part of a set of nine concentrated solar power simulation problems~\cite{MScMLG} that serve as a benchmark for blackbox optimization solvers, available at \linebreak \href{https://github.com/bbopt/solar}{\tt github.com/bbopt/solar}.
The \solar{1} problem aims at maximizing the heliostat field energy output throughout one day under constraints of field geometry and cost. The problem is noisy and has several local optima. It has $n=9$ variables, among which one is discrete, and $m=5$ constraints. The \solar{1} problem has the specificity of having two additional adjustable parameters: the seed and the number of replications. Since the problem contains some stochasticity, the seed for the random generator can be chosen. It is \textit{not} the same seed as the one of \nomad. Two instances of \solar{1} with two different seeds are considered as two different problems in this work. The other parameter enables to replicate the evaluations in order to smooth the problem and compensate the noise. It was fixed to ten in the experiments.

In order to generate several problem instances, fifteen different seeds are chosen for the problem - not for \nomad - instead of multiple starting points. Because the problem is especially expensive, \nomad has been tested with only one seed, and not four, in order to reduce the number of optimization runs. The evaluation budget is $800(n+1)$. As with the previous problem, only the best formulations from the aircraft range problem have been tested. The best one among them is SP8 for both smooth and nonsmooth alternatives. They manage to yield better results than \textsf{no search}, however, \textsf{quad search} is clearly the best algorithm on this problem.

\begin{figure}[htb]
\begin{subfigure}{.5\textwidth}
  \centering
  \includegraphics[width=\linewidth]{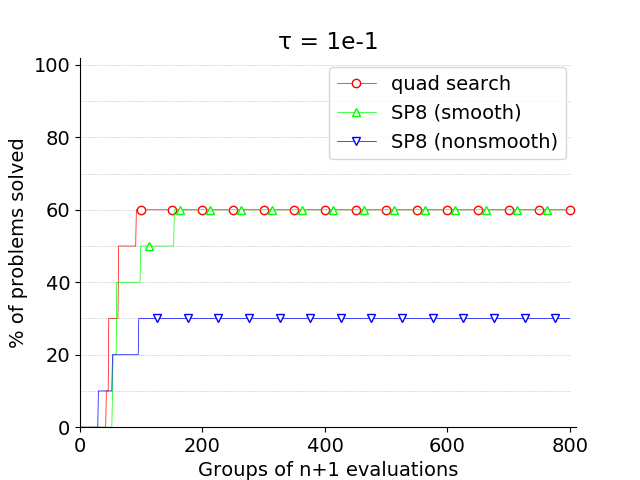}
\end{subfigure}
\begin{subfigure}{.5\textwidth}
  \centering
  \includegraphics[width=\linewidth]{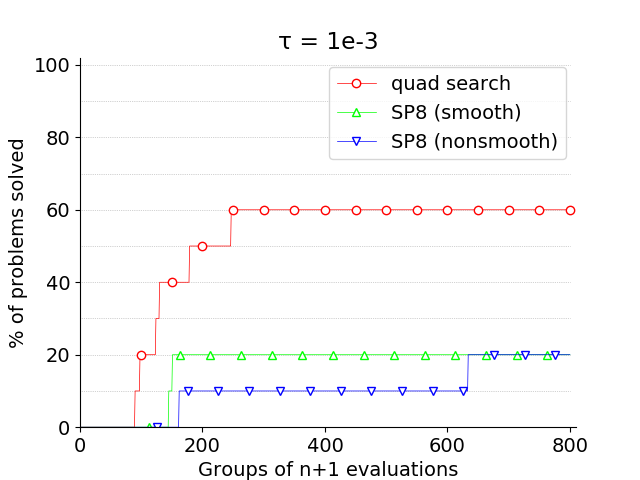}
\end{subfigure}
\begin{subfigure}{.5\textwidth}
  \centering
  \includegraphics[width=\linewidth]{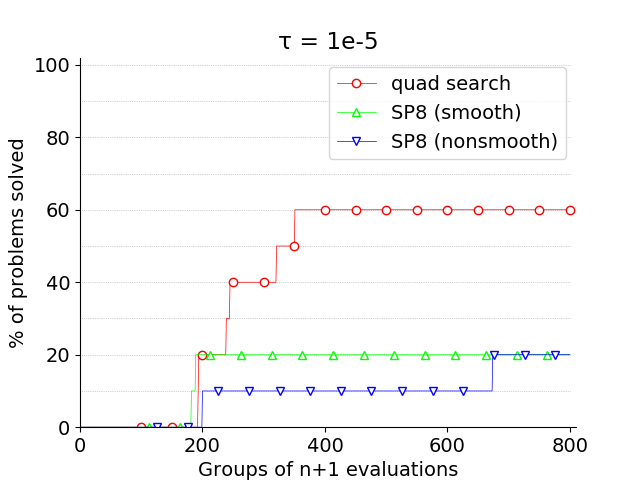}
\end{subfigure}
\begin{subfigure}{.5\textwidth}
  \centering
  \includegraphics[width=\linewidth]{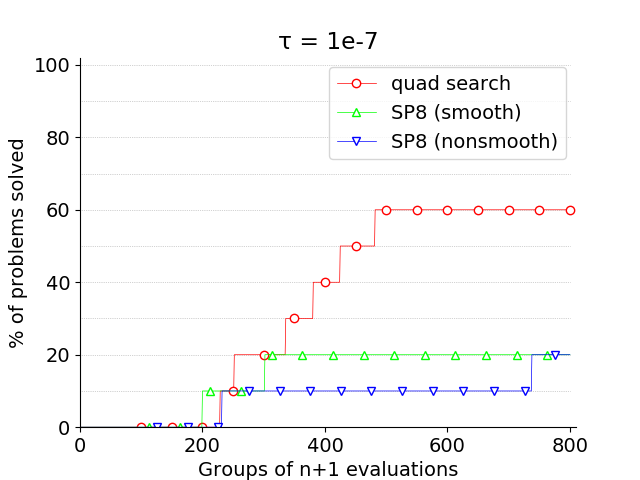}
\end{subfigure}
\caption{Data profiles. \textsf{quad search} vs. SP8 with smooth uncertainty vs. SP8 with nonsmooth uncertainty on \solar{1}.}
\label{fig:quadsearch_vs_SP8smooth_vs_SP8nonsmooth}
\end{figure}

Figure~\ref{fig:quadsearch_vs_SP8smooth_vs_SP8nonsmooth} shows the data profiles of SP8 with smooth and nonsmooth alternatives along with \textsf{quad search}. \dfn and \textsf{no search} are not represent because they do not manage to solve one problem even at the largest tolerance. At tolerance $\tau=10^{-1}$, the smooth alternative manages to solve as many problems as \textsf{quad search}. However, for higher tolerances, the latter is by far the best alternative. The poor performance of the extended aggregate models can be attributed to the stochasticity of the problem that importantly deteriorates the models. Nonetheless, it also impacts \textsf{quad search}. The superior performance of the latter is also due to the nature of the constraints: three out of five are linear, and one is cubic, which gives a non negligible advantage to quadratic models.

%-----------------------------------%
\subsection[The styrene problem]{The \styrene problem}
%-----------------------------------%

The \styrene problem is a chemical engineering simulator for styrene production described in~\cite{AuBeLe08}
and available at \href{https://github.com/bbopt/styrene}{\tt github.com/bbopt/styrene}. It aims at maximizing the net present value of the styrene production process under structural, chemical and financial constraints with variables comprised of physical parameters and structure. The problem is deterministic but nonsmooth and with omnipresent hidden constraints, i.e., the simulation often fails to return a value even when all constraints are met. A random sampling resulted in almost 60\% of failures in~\cite{GrLeD2011}. In addition, four constraints are binary. The problem has $n=8$ bounded variables and $m=11$ constraints.

Feasible regions may by especially hard to find on this problem. Consequently, twelve starting points were generated in three relatively easy regions. The evaluation budget is $600(n+1)$. Like the previous problem \solar{1}, the same subset of the most promising formulations found on the aircraft range problem has been tested. The best formulations are SP$1_{1}$ with the smooth alternative and SP$3_{0.01}$ with the nonsmooth alternative.

\begin{figure}[!ht]
\begin{subfigure}{.5\textwidth}
  \centering
  \includegraphics[width=\linewidth]{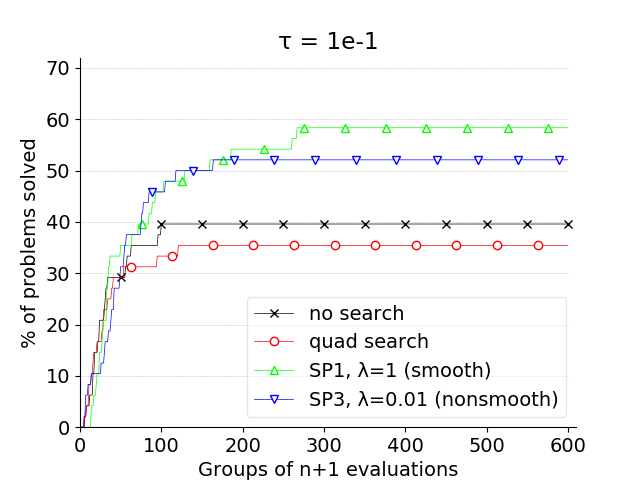}
\end{subfigure}
\begin{subfigure}{.5\textwidth}
  \centering
  \includegraphics[width=\linewidth]{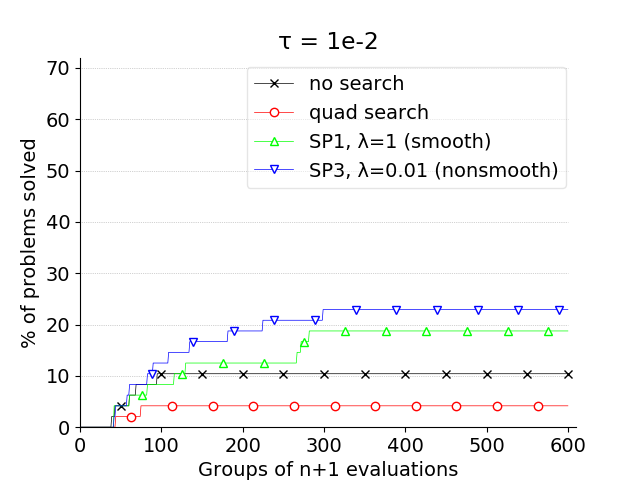}
\end{subfigure}
\caption{Data profiles. \textsf{no search} vs. \textsf{quad search} vs. SP$1_{1}$ with smooth uncertainty vs. SP$3_{0.01}$ with nonsmooth uncertainty on \styrene.}
\label{fig:nosearch_vs_quadsearch_vs_SPsmooth1_vs_SP3nonsmooth0.01}
\end{figure}

Figure~\ref{fig:nosearch_vs_quadsearch_vs_SPsmooth1_vs_SP3nonsmooth0.01} shows the data profiles of \textsf{no search}, \textsf{quad search}, SP$1_{1}$ with the smooth uncertainty and SP$3_{0.01}$ with the nonsmooth uncertainty at tolerance $\tau=10^{-1}$ and $\tau=10^{-2}$. On this problem, the extended aggregate models perform better than both \textsf{no search} and \textsf{quad search}. Besides, unlike all the others problems the latter yields worse results than \textsf{no search} due to the binary constraints. \dfn struggles to find feasible solutions and the results are not presented.
%

%-----------------------------------%
\subsection{Results of \shebo}
\label{subsec:results_shebo}
%-----------------------------------%

Unlike the other algorithms in this work, \shebo does not take into account a starting point as an input, thereby making the comparison through data profiles impossible. In order to analyze the results, the best value found for a given optimization run is denoted by $f^*$, the real time (in minutes) needed to reach $f^*$ is denoted by $t^*$, and the number of evaluations needed to reach $f^*$ is denoted by $k^*$. Table~\ref{tab:results_shebo} shows for most of the solvers seen before the median $f^*$, $t^*$ and $k^*$ on each problem considering all starting points and seeds, denoted by $f^*_{\mathrm{m}}$, $t^*_{\mathrm{m}}$ and $k^*_{\mathrm{m}}$, respectively, as well as the best $f^*$ found on all runs, denoted by $f^*_{\mathrm{best}}$, and the total real time needed for all runs, denoted by $t_{\mathrm{tot}}$. For every problem, except the analytical problems, the three quantities are shown for \textsf{no search}, \textsf{quad search}, the best formulation found with the smooth uncertainty, the best formulation found with the nonsmooth uncertainty, \dfn and \shebo. However, \shebo has been run only once on each problem due to its long running time. Consequently, all the values shown on the lines corresponding to \shebo relate only to a single run, and $f^*_{\mathrm{m}}$ is equal to $f^*_{\mathrm{best}}$.

\begin{table}
\small
\centering
\renewcommand{\arraystretch}{1.5}
\renewcommand{\tabcolsep}{1.5pt}
\begin{tabular}{lllllllllll} 
\toprule
\multicolumn{1}{l|}{} & \multicolumn{5}{c|}{\begin{tabular}[c]{@{}c@{}}Aircraft range\\(40 runs except for \shebo)\end{tabular}} & \multicolumn{5}{c|}{\begin{tabular}[c]{@{}c@{}}Simplified wing\\(40 runs except for \shebo)\end{tabular}} \\ 
\cline{2-11}
\multicolumn{1}{l|}{} & $f^*_{\mathrm{m}}$ & $t^*_{\mathrm{m}}$ & $k^*_{\mathrm{m}}$ & $f^*_{\mathrm{best}}$ & \multicolumn{1}{l|}{$t_{\mathrm{tot}}$} & $f^*_{\mathrm{m}}$ & $t^*_{\mathrm{m}}$ & $k^*_{\mathrm{m}}$ & $f^*_{\mathrm{best}}$ & \multicolumn{1}{l|}{$t_{\mathrm{tot}}$} \\ 
\hline
\textsf{no search}   & -3964.204696 & 1   & 4140 & -3964.204700 & 118  & -16.4059 & 3  & 2405 & -16.6119 & 397 \\
\textsf{quad search} & -3964.204698 & 5   & 3064 & -3964.204701 & 834  & -16.5808 & 12 & 1291 & -16.6120 & 2181 \\
smooth               & -3964.204699 & 10  & 3218 & -3964.204701 & 1836 & -16.6063 & 4  & 2115 & -16.6120 & 921 \\
nonsmooth            & -3964.204699 & 9   & 2555 & -3964.204700 & 1614 & -16.5924 & 2  & 836  & -16.6120 & 393 \\
\dfn                 & -1749.256329 & 0.1 & 73   & -3143.101404 & 18   & -15.0774 & 3  & 1134 & -16.6112 & 588 \\
\shebo               & $\quad\quad\quad\quad\cdot$ & 784 & 4404 & -3723.074357 & 3210 & $\quad\quad\quad\cdot$ & 64 & 2190 & -16.5501 & 356 \\
\cline{2-11}
\multicolumn{1}{l|}{} & \multicolumn{5}{c|}{\begin{tabular}[c]{@{}c@{}}\solar{1}\\(10 runs except for \shebo)\end{tabular}} & \multicolumn{5}{c|}{\begin{tabular}[c]{@{}c@{}}\styrene\\(48 runs except for \shebo)\end{tabular}} \\ 
\cline{2-11}
\multicolumn{1}{l|}{} & $f^*_{\mathrm{m}}$ & $t^*_{\mathrm{m}}$ & $k^*_{\mathrm{m}}$ & $f^*_{\mathrm{best}}$ & \multicolumn{1}{l|}{$t_{\mathrm{tot}}$} & $f^*_{\mathrm{m}}$ & $t^*_{\mathrm{m}}$ & $k^*_{\mathrm{m}}$ & $f^*_{\mathrm{best}}$ & \multicolumn{1}{l|}{$t_{\mathrm{tot}}$} \\ 
\hline
\textsf{no search}   & -508184.0 & 190 & 4064 & -660888.3 & 2384 & -29305150 & 12  & 1541 & -33613200 & 3611 \\
\textsf{quad search} & -723623.2 & 510 & 3795 & -835124.4 & 6837 & -29301200 & 41  & 1675 & -33000800 & 9750 \\
smooth               & -677669.7 & 305 & 3688 & -849192.5 & 6303 & -32704300 & 49  & 1809 & -33705600 & 19015 \\
nonsmooth            & -658837.2 & 224 & 3359 & -805799.5 & 6190 & -32235250 & 48  & 1647 & -33697400 & 16086 \\
\dfn                 & -332582.8 & 68  & 1448 & -391849.8 & 735  & -22851550 & 0.3 & 38   & -28517600 & 49 \\
\shebo               & $\quad\quad\quad\cdot$ & 2865 & 6190 & -815086.2 & 4475 & $\quad\quad\quad\cdot$ & 741 & 5365 & -32873400 & 774 \\
\bottomrule
\end{tabular}
\caption{Results of \shebo compared to other algorithms.}
\label{tab:results_shebo}
\end{table}

For \solar{1} and \styrene, \shebo has a better $f^*_{\mathrm{best}}$ than the median $f^*$ of all algorithms. However, on \solar{1} the proposed approach, both smooth and nonsmooth, manages to find a better $f^*_{\mathrm{best}}$ over all 10 runs after only 1.5 times more real time than the single run of \shebo. On \styrene, all the \mads algorithms manage to find a better $f^*_{\mathrm{best}}$ than \shebo after all 48 runs, with significantly more real time though. On the aircraft range and the simplified wing problems, the proposed approach not only finds a better $f^*_{\mathrm{best}}$ than \shebo but also a better $f^*_{\mathrm{med}}$ after only a fraction of the time per run. On the simplified wing problem, \dfn also finds a better $f^*_{\mathrm{best}}$ than \shebo with 1.7 times more real time. Regarding the number of evaluations, \shebo requires more function evaluations to reach its best value than the other algorithms, except on the simplified wing problem. Overall, on the present experiments, one single run with \shebo guaranties a decent value, only it demands a much larger computation time which is better invested running the proposed approach multiple times.

One final remark can be made about the real time required by the extended aggregate models compared to the quadratic models. The proposed approach is not always longer than \textsf{quad search}, both in terms of median time and total time. The quadratic models indeed tend to result in very long runs when they find a good bassin of solutions. Conversely, they result in very short runs when then they do not manage to find a good solution, depending on the starting point.

%-------------------------------------------%
\section{Discussion}
\label{sec:conclusion}
%-------------------------------------------%

This work proposes an extension to ensembles of models that enables to compute an uncertainty at any given point. The resulting extended aggregate models behave like stochastic models, i.e., they produce at any given point $x$ a prediction $\hat f(x)$ and also an uncertainty $\hat\sigma(x)$, thus enabling to use tools inspired by Bayesian optimization. The proposed extended aggregate models are incorporated into the search step of \mads where at each iteration a surrogate subproblem derived from Bayesian optimization is solved in order to come up with new candidate points. The proposed approach may be used in any direct search method based on the search-poll paradigm, or in any approach akin to efficient global optimization if adapted. Any ensemble of models can be used along with any weight attribution technique provided that at least two models have a strictly positive weight at any moment.

The resulting algorithm has been tested on seven analytical problems, two multidisciplinary optimization problems and two simulation problems. The results show that the proposed extended aggregate models incorporated into \mads find better solutions than \mads without search step or with the help of quadratic models on three expensive problems out of four. They also find better solutions than the stochastic models that they replace while requiring much less computational time. It should be noted that the models used to build the aggregate models must remain moderately expensive to compute, otherwise the method might loose its advantage in terms of computation time. The proposed approach does not show an advantage over quadratic models on analytical problems. In addition, the latter yield better results on the \solar{1} problem which most of the constraints are linear. The comparison to other solvers shows that the proposed approach has a clear advantage over \dfn, and is more interesting than \shebo when given the same computation time. An extended study of the various sub-problem formulations has not been conducted but based on the present results the formulations SP1, SP3 and SP8 can be recommended.

Future work may explore the uncertainty for the constraint independently from that of the objective, e.g., smooth uncertainty for the objective together with nonsmooth uncertainty for the constraint, since in this work the two were coupled. Other weight attribution techniques than that described in 
Section~\ref{subsec:weight_attribution} may also be considered. The influence of the simplex used to build simplex gradients in~\eqref{eq:uncertainty_obj_smooth} as well as the positive spanning set in~\eqref{eq:uncertainty_obj_nonsmooth} have not been studied in this work. Besides, the parameter $\lambda$ of the surrogate subproblems has been carefully selected for each formulation but remains constant over the optimization once determined. It might instead be dynamically updated depending on the result of the search or merely follow a predetermined trend like decreasing with the number of iterations. The parameter $\alpha$ in~\eqref{eq:uncertainty} is proportional to the global variance of the cache. It could be refined in order to represent local trends better, for instance by taking into account the values of the cache only in a restricted area around the evaluated point, or by removing outliers. Finally, the formulations were chosen on a purely empirical basis. Little effort has been made to finely understand the behaviour and the performance thereof. A thorough analysis of the benefits of each formulations in the context of extended aggregate models might be a judicious undertaking.

%-------------------------------------------%
%-------------------------------------------%
\bibliographystyle{plain}
\bibliography{bibliography}
%-------------------------------------------%

\appendix

%-------------------------------------------%
\section{Positive spanning set and simplex construction}\label{appendix:pps_simplex}
%-------------------------------------------%

The simplex and the positive spanning set described below are built in the \textit{scaled} search space. Before constructing the models, the \nomad software used in this work scales each input variable $x_i$, $i\in\{1,2,\dots,n\}$, using the mean and the variance of the points of the cache. This is done to give the same importance to all the variables regardless of their initial amplitude. Consequently, the simplex and the positive spanning set are isotropic in the scaled search space, but not in the actual search space.

\begin{itemize}

\item
The simplex centred on $x\in\mathbb{R}^n$ used to build the simplex gradients $\nabla_S f(x)$ in Equation~\eqref{eq:uncertainty_obj_smooth} is
$$
\left\{ x + 0.001 d_i\ :\ i\in\{1,2,\dots,n+1\} \right\}
$$
where the directions $d_i$, $i\in\{1,2,\dots,n+1\}$, are constructed according to the following procedure:
$$
\begin{aligned}
d_i =\left\{\begin{array}{ll}
     e_i - \displaystyle\frac{1+\frac{1}{\sqrt{n+1}}}{n}\times [1,1,\dots,1]^{\top}, & \mbox{if } i\in\{1,2,\dots,n\}   \\
     \displaystyle\frac{1}{\sqrt{2(n+1)}} \times  [1,1,\dots,1]^{\top}, & \mbox{if } i=n+1
\end{array} \right. 
\end{aligned}
$$
where $e_i$ is the vector $[0,\dots,1,\dots,0]^{\top}$ with value~1 at the $i$th position. This procedure forms a regular simplex centred on $x$ with side length equal to $\sqrt{2}$ in any dimension. The factor 0.001 was chosen empirically on preliminary tests.

\item
The positive spanning set centred on $x\in\mathbb{R}^n$ used in
Equation~\eqref{eq:uncertainty_obj_nonsmooth} is
$$
\left\{ x\pm 0.005 e_i\ :\ i\in\{1,2,\dots,n\} \right\}
$$
This positive spanning set contains $2n$ elements. The factor 0.005 was chosen empirically on preliminary tests.

\end{itemize}

%--------------------------------------------%
\section{Surrogate subproblem formulations}
\label{appendix:SP_formulations}
%--------------------------------------------%

\begin{align*}
\tag{SP1-F$\sigma$}
    \min_{x\in\mathcal{X}}&\ \ \hat f(x)-\lambda\hat\sigma_f(x) \\
    \mathrm{s.t.}&\ \ \hat c_j(x)-\lambda\hat\sigma_j(x)\leq0,\ \ j=1,2,\dots,m \\
    \\
\tag{SP2-F$\sigma$P}
    \min_{x\in\mathcal{X}}&\ \ \hat f(x)-\lambda\hat\sigma_f(x) \\
    \mathrm{s.t.}&\ \ \mathrm{P}(x)\geq p_c \\
    \\
\tag{SP3-EI$\sigma$}
    \min_{x\in\mathcal{X}}&\ -\mathrm{EI}(x)-\lambda\hat\sigma_f(x) \\
    \mathrm{s.t.}&\ \ \hat c_j(x)-\lambda\hat\sigma_j(x)\leq0,\ \ j=1,2,\dots,m \\
    \\
\tag{SP4-EFI}
    \min_{x\in\mathcal{X}}&\ -\mathrm{EFI}(x) \\
    \\
\tag{SP5-EFI$\sigma$}
    \min_{x\in\mathcal{X}}&\ -\mathrm{EFI}(x)-\lambda\hat\sigma_f(x) \\
    \\
\tag{SP6-EFI$\mu$}
    \min_{x\in\mathcal{X}}&\ -\mathrm{EFI}(x)-\lambda\hat\sigma_f(x)\mu(x) \\
    \\
\tag{SP7-EFIC$\mu$}
    \min_{x\in\mathcal{X}}&\ -\mathrm{EFI}(x)-\lambda(\mathrm{EI}(x)\mu(x)
    +\mathrm{P}(x)\hat\sigma_f(x)) \\
    \\
\tag{SP8-PFI}
    \min_{x\in\mathcal{X}}&\ -\mathrm{PFI}(x)
\end{align*}

\end{document}